\documentclass{amsart}
\usepackage{amsfonts,amssymb,amscd,amsmath,enumerate,verbatim,xypic}
\xyoption{curve}
\usepackage{calc}
\CompileMatrices

\newcommand{\ts}{\textstyle}
\newcommand{\sst}{\scriptstyle}
\newcommand{\sss}{\scriptscriptstyle}

\newcommand{\ges}{{\sss\geqslant}}
\newcommand{\les}{{\sss\leqslant}}

\newcommand{\sles}{{\sss<}}

\newcommand{\aq}[4]{\operatorname{D}_{#1}(#2\hskip.7pt|\hskip.7pt#3;#4){}}
\newcommand{\aqdim}[2]{\operatorname{AQ-dim}_{#2}{#1}}

\newcommand{\barr}[2]{\operatorname{B}_{#1}(#2)}
\newcommand{\bark}[2]{\ov{\operatorname B}_{#1}(#2)}
\newcommand{\hark}[2]{\operatorname{H}\!\ov{\operatorname B}_{#1}(#2)}
\newcommand{\bsa}{{\boldsymbol a}}
\newcommand{\bsf}{{\boldsymbol f}}
\newcommand{\bsg}{{\boldsymbol g}}

\newcommand{\bsx}{{\boldsymbol x}}
\newcommand{\bu}{{\sss\bullet}}
\newcommand{\card}{\operatorname{card}}
\newcommand{\chr}{\operatorname{char}}
\newcommand{\col}{\colon}
\newcommand{\dd}{\partial}
\newcommand{\depth}{\operatorname{depth}}
\newcommand{\ED}[3]{{}^{#1}\!{d}_{#2,#3}}
\newcommand{\EH}[3]{{}^{#1}\!\operatorname{E}_{#2,#3}}
\newcommand{\FF}{\operatorname{F}}
\newcommand{\dvip}[2]{{{#2}^{(#1)}}}
\newcommand{\edim}{\operatorname{edim}}
\newcommand{\fd}{\operatorname{fd}}
\newcommand{\fa}{{\mathfrak a}}
\newcommand{\fb}{{\mathfrak b}}
\newcommand{\fm}{{\mathfrak m}}
\newcommand{\fn}{{\mathfrak n}}
\newcommand{\fp}{{\mathfrak p}}
\newcommand{\fq}{{\mathfrak q}}
\newcommand{\gr}{\operatorname{gr}}
\newcommand{\G}{$\Gamma$}
\newcommand{\vG}{$\varGamma$}
\newcommand{\gin}[1]{{\Gamma}\!\!\operatorname{-ind}(#1)}
\newcommand{\Gin}[2]{{\Gamma}\!\!\operatorname{-ind}_{#1}(#2)}
\newcommand{\hch}{\operatorname{HH}}
\newcommand{\hh}[1]{\operatorname{H}(#1)}

\newcommand{\HH}[2]{\operatorname{H}_{#1}(#2)}
\newcommand{\hra}{\hookrightarrow}
\newcommand{\id}{\operatorname{id}}
\newcommand{\im}{\operatorname{Im}}
\newcommand{\ind}[1]{{\operatorname{ind}(#1)}}
\newcommand{\Ind}[2]{{\operatorname{ind}_{#1}(#2)}}
\newcommand{\Ker}{\operatorname{Ker}}
\newcommand{\la}{\operatorname{\!\langle\!}}
\newcommand{\length}{\operatorname{length}}

\newcommand{\lra}{\longrightarrow}
\newcommand{\Max}{\operatorname{Max}}
\newcommand{\nat}{{}^\natural}
\newcommand{\ov}{\overline}
\newcommand{\pie}[2]{{\pi_{#1}(#2)}}
\newcommand{\pies}[1]{{\pi_{\bu}(#1)}}
\newcommand{\polreg}{\operatorname{pol\,reg}}
\newcommand{\rank}{\operatorname{rank}}
\newcommand{\ra}{\operatorname{\!\rangle}}
\newcommand{\sgn}[2]{{(-1)^{|#1||#2|}}}
\newcommand{\Shift}{{\Sigma}}
\newcommand{\Spec}{\operatorname{Spec}}
\newcommand{\susp}{{\sst\Sigma}}
\newcommand{\sym}{\operatorname{Sym}}
\newcommand{\Tor}[4]{\operatorname{Tor}_{#1}^{#2}(#3,#4){}}
\newcommand{\rTor}[4]{{\operatorname{tor}}{}_{#1}^{#2}(#3,#4){}}
\newcommand{\tra}{\twoheadrightarrow}
\newcommand{\rat}{\rightarrowtail}
\newcommand{\ul}{\underline}
\newcommand{\var}{{\hskip.7pt\vert\hskip.7pt}}
\newcommand{\ve}{{\varepsilon}}
\newcommand{\vf}{{\varphi}}
\newcommand{\wcat}{\operatorname{wcat}}
\newcommand{\wh}{\widehat}
\newcommand{\wt}{\widetilde}

\newcommand{\xra}{\xrightarrow}
\newcommand{\ZZ}{\operatorname{Z}}

\newcommand{\BN}{{\mathbb N}}
\newcommand{\BQ}{{\mathbb Q}}
\newcommand{\BZ}{{\mathbb Z}}

\swapnumbers \theoremstyle{plain}
\newtheorem{theorem}{Theorem}[section]
\newtheorem{proposition}[theorem]{Proposition}
\newtheorem{lemma}[theorem]{Lemma}
\newtheorem{corollary}[theorem]{Corollary}
\newtheorem*{main}{Theorem I}
\newtheorem*{second}{Theorem II}

\theoremstyle{definition}

\newtheorem{definition}[theorem]{Definition}
\newtheorem{example}[theorem]{Example}
\newtheorem*{exa}{Examples}
\newtheorem{chunk}[theorem]{}
\newtheorem{remark}[theorem]{Remark}

\swapnumbers \theoremstyle{plain}

\theoremstyle{remark}
\newtheorem*{conjecture}{Quillen's Conjecture}

\begin{document}

{\ }
\vskip-2.4truecm
\noindent{\small{\sc Ann. Sci. \'Ecole Norm. Sup. (4)}

\noindent{(to appear)}}
\bigskip
\bigskip
\bigskip
\medskip

\title[Homology of algebra retracts]
{Andr\'e-Quillen homology of algebra retracts}

\date{October 14, 2001; revised April 14, 2002} %,\ \printtime}

\author[L.~L.~Avramov]{Luchezar L.~Avramov}
\address{Department of Mathematics,
Purdue University, West Lafayette, IN~47907, USA}
\email{avramov@math.unl.edu}
\curraddr{Department of Mathematics, University of Nebraska,
Lincoln, NE 68588, USA}
\author[S.~Iyengar]{Srikanth Iyengar}
\address{Pure Mathematics, Hicks Building,
University of Sheffield, Sheffield S3 7RH, UK}
\email{iyengar@math.missouri.edu}
\curraddr{Department of Mathematics, University of Missouri,
Columbia, MO 65201, USA}
\thanks{L.L.A. was partly supported by a grant from the N.S.F.
\endgraf
S.I. was supported by a grant from the E.P.S.R.C}
  \maketitle

\noindent{\bf Abstract.}
Given a homomorphism of commutative noetherian rings $\vf\col R\to S$,
Daniel Quillen conjectured in 1970 that if the Andr\'e-Quillen
homology functors $\aq nSR-$ vanish for all $n\gg0$, then they vanish for
all $n\ge3$.  We prove the conjecture under the additional hypothesis
that there exists a homomorphism of rings $\psi\col S\to R$ such that
$\vf\circ\psi=\id_S$.  More precisely, in this case we show that $\psi$
is complete intersection at $\vf^{-1}(\fn)$ for every prime ideal $\fn$
of $S$.  Using these results, we describe all algebra retracts $S\to
R\to S$ for which the algebra $\Tor\bu RSS$ is finitely generated over
$\Tor0RSS=S$.

\bigskip
\bigskip
\bigskip
\bigskip
\bigskip
\bigskip

\centerline
{\bf HOMOLOGIE D'ANDR\'E-QUILLEN DES ALG\`EBRES SCIND\'EES}

\bigskip
\bigskip
\bigskip

\noindent{\bf Resum\'e.}
\'Etant donn\'e un homomorphisme $\vf\col R\to S$ d'anneaux commutatifs
noeth\'eriens, Daniel Quillen a conjectur\'e en 1970 que si les foncteurs
$\aq nSR-$ d'homologie d'Andr\'e-Quillen sont nuls pour tout $n\gg0$,
alors ils sont nuls pour tout $n\ge3$.  Nous d\'emontrons cette conjecture
sous l'hypoth\`ese suppl\'ementaire qu'il existe un homomorphisme
d'anneaux $\psi\col S\to R$ tel que $\vf\circ\psi=\id_S$.  Plus
pr\'ecisemment, nous montrons que dans ce cas $\psi$ est d'intersection
compl\`ete en $\vf^{-1}(\fn)$ pour tout id\'eal premier $\fn$ de $S$.
En utilisant ces r\'esultats, nous d\'ecrivons toutes les alg\`ebres
scind\'ees $S\to R\to S$ pour lesquelles l'alg\`ebre $\Tor\bu RSS$
est finiment engendr\'ee sur $\Tor0RSS=S$.

\newpage

\section*{Introduction}
\setcounter{section}{0}

Let $\vf\col R\to S$ be a homomorphism of commutative noetherian rings.

For each $n\ge0$, let $\aq nSR-$ denote the $n$th cotangent homology
functor on the category of $S$-modules, defined by Andr\'e \cite{simp}
and Quillen \cite{Qu}. To study how vanishing of these Andr\'e-Quillen
homology functors relates to the structure of $\vf$, we define the {\em
Andr\'e-Quillen dimension\/} of $S$ over $R$ to be the number
\[
\aqdim SR = \sup\{n\in\BN\mid \aq nSR-\ne 0\}\,;
\]
in particular, $\aqdim SR=-\infty$ if and only if $\aq nSR-=0$ for all
$n\in\BZ$.

Vanishing of Andr\'e-Quillen homology in low dimensions characterizes
important classes of homomorphisms of noetherian rings. Recall that $\vf$
is {\em regular\/} if it is flat with geometrically regular fibers.  It is
{\em \'etale\/} if, in addition, it is of finite type and unramified.
A general {\em locally complete intersection\/}, or {\em l.c.i.\/},
property is defined in \ref{def: lci}; when $\vf$ is of finite type,
it means that in some (equivalently, every) factorization of $\vf$ as an
inclusion into a polynomial ring followed by a surjection, the kernel of
the second map is locally generated by a regular sequence.  The following
results were proved in \cite{simp}, \cite{Qu} for maps $\vf$ of finite
type, and in \cite{lis}, \cite{lci} in general:
\begin{enumerate}[{\rm(A)}]
\item
$\aqdim SR=-\infty$ and $\vf$ is of finite type if and only if
$\vf$ is \'etale.
\item
$\aqdim SR \leq 0$ if and only if $\aq 1SR-=0$, if and only if
$\vf$ is regular.
\item
$\aqdim SR \leq 1$ if and only if $\aq 2SR-=0$, if and only $\vf$ is
l.c.i.
 \end{enumerate}

Further research on homomorphisms of finite Andr\'e-Quillen dimension has
been driven by two conjectures, stated by Quillen in 1970.  One of them,
\cite[(5.7)]{Qu}, is for maps {\em locally of finite flat dimension\/}:
For each prime ideal $\fn$ of $S$ the $R$-module $S_\fn$ has a finite
resolution by flat $R$-modules.  That conjecture was proved in \cite{lci}:
\begin{enumerate}[{\rm(A)}]
\item[{\rm(D)}]
$\aqdim SR<\infty$ and $\vf$ is locally of finite flat dimension
if and only if $\vf$ is l.c.i.
 \end{enumerate}

As a consequence, if $\vf$ is locally has finite flat dimension, then
$\aqdim SR <\infty$ implies $\aqdim SR\le1$.  The remaining conjecture,
\cite[(5.6)]{Qu}, predicts the behavior of Andr\'e-Quillen dimension
when no flatness hypothesis is available.

\begin{conjecture}
If $\aqdim SR<\infty$, then $\aqdim SR\leq 2$.
 \end{conjecture}

No structure theorem is known for $R$-algebras $S$ with $\aqdim SR\leq 2$,
so the conjecture presents a significant challenge beyond the generic
difficulty of computing the modules $\aq nSRM$, defined in terms of
simplicial resolutions.  This partly explains why so few cases have
been settled.  In \cite{lci} the conjecture is proved when one of the
rings $R$ or $S$ is locally complete intersection.  Indirect evidence is
obtained in \cite{acy}: If $\vf$ is a large homomorphism of local rings
in the sense of \cite{Le2}, $R$ has characteristic $0$, and $\aqdim SR$
is an odd integer, then $\aqdim SR=1$.

Our main result establishes Quillen's Conjecture when $S$ is an {\em
algebra retract\/} of $R$, meaning that there exists a homomorphism of
rings $\psi\col S\to R$ such that $\vf\circ\psi=\id_S$; any homomorphism
$\psi$ with this property is called a {\em section\/} of $\vf$.  Algebra
retracts frequently arise from geometric considerations.  For instance, to
study a morphism of schemes $X\to Y$ one often uses the induced diagonal
embedding $X\to X\times_YX$.  The underlying algebraic construction is the
homomorphism of rings $\vf\col S\otimes_AS\to S$ defined by $\vf(s'\otimes
s'')=s's''$; the ring $S$ is an algebra retract of $R=S\otimes_AS$,
with section $\psi(s)=s\otimes1$. A different type of retracts arises in
constructions of projective schemes.  They typically involve a graded
$S$-algebra $R=\bigoplus_{i=0}^\infty R_i$ with $R_0=S$; the relevant
homomorphisms $\vf$ and $\psi$ are, respectively, the canonical surjection
$R\to (R/R_{\ges 1})=S$ and the inclusion $S=R_0\subseteq R$.

An important aspect of our result is that it connects the homological
conditions in the conjecture through the structure of retracts
of finite Andr\'e-Quillen dimension.  Let, as always, $\Spec S $ denote
the set of prime ideals of $S$.  If $\vf$ has a section $\psi$, then
for every $\fn\in\Spec S$ one can find a set $\bsx$ of formal
indeterminates over $S_\fn$ and an ideal $\fb$ contained in $\fn(\bsx)+
(\bsx)^2$ that fit into a commutative diagram
\begin{equation*}
\label{factoriz}
\xymatrixrowsep{2pc}
\xymatrixcolsep{2.5pc}
\begin{gathered}
\xymatrix{
&{\ }(R_{\fm})^*
\ar@{->}[dr]^{(\varphi_\fn)^*}
\\
S_\fn
\ar@{->}[r]^{\psi'}
\ar@{->}[ur]^{(\psi_\fm)^*}
& {S_\fn[[\bsx]]}/{(\fb)}
\ar@{->}[u]^{\cong}
\ar@{->}^{\vf'}[r]
& S_\fn
}
\end{gathered}
\tag{E${}_\fn$}
\end{equation*}
of homomorphisms of rings, where $\fm=\vf^{-1}(\fn)$, asterisks ${\,}^*$
denote $(\Ker(\vf))$-adic completion, $\psi'$ is the natural injection
and $\vf'$ the surjection with kernel $(\bsx)$.

For every real number $c$ set $\lfloor c\rfloor=\sup\{i\in\BZ\mid
i\le c\}$.

\begin{main}
Let $\vf\col R\to S$ be a homomorphism of rings and set $\fa=\Ker(\vf)$.
If $\vf$ admits a section and $R$ is noetherian, then the following
conditions are equivalent.
\begin{enumerate}[{\quad\rm(i)}]
\item
$\aqdim SR<\infty$\,.
\item
$\aqdim SR\leq 2$\,.
\item
$\aq 3SR-=0$\,.
\item
$\aq nSR-=0$ for some $n\ge 3$ such that $\lfloor \frac {n-1}2 \rfloor!$
is invertible in $S$\,.
\item
For each $\fn\in\Spec S$, the ideal $\fb$ in some (respectively, every)
commutative diagram \eqref{factoriz} is generated by a regular sequence.
 \end{enumerate}
  \end{main}

We apply the results discussed above in concrete cases, illustrating the
known fact that all dimensions allowed under Quillen's Conjecture do occur.

\begin{exa}
Let $x$, $y$ be indeterminates over $S$.  The natural homomorphisms
\[
\xymatrixrowsep{.5pc}
\xymatrixcolsep{1.5pc}
\xymatrix{
S
\ar@{->}[r]
&S[x,y]
\ar@{->}[r]\ar@{=}[d]
&S[x,y]/(x^2,xy,y^2)
\ar@{->}[r]\ar@{=}[d]
&S[x]/(x^2)
\ar@{->}[r]\ar@{=}[d]
&S
\\
&P
&R
&T
}
\]
provide the following list of Andr\'e-Quillen dimensions:
\begin{alignat*}{4}
\aqdim SS&=\aqdim PP&&=&&\aqdim RR&&=\aqdim TT=-\infty\\
\aqdim PS&=0\\
\aqdim TS&=\aqdim TP&&=&&\aqdim SP&&=\,1 \\
\aqdim ST&=2\\
\aqdim RS&=\aqdim RP&&=&&\aqdim TR&&=\aqdim SR =\infty
\end{alignat*}
Indeed, (A), (B), and (C) yield the equalities in the first three lines;
(C) also implies $\aqdim ST\ge2$.  Because $S$ is a retract of $T$,
Theorem I provides the converse inequality; since $T$ and $S$ are
retracts of $R$, the theorem also computes the last two dimensions on
the last line.  The two remaining dimensions on that line are given by
(D), because $R$ has finite flat dimension over $S$ and over $P$.
 \end{exa}

We use Theorem I together with our results in \cite{hh} in a situation
that does not {\em a priori\/} involve Andr\'e-Quillen homology---the
classical homology of an algebra retract $S\to R\to S$.  In that case
$\Tor0RSS=S$ and $\Tor\bu RSS$ is a graded-commutative algebra with
divided powers, but precise information on its structure is available in
two instances only: when $S$ is a field, cf.\ \cite{Le}, \cite{GL},
or when $R\to S$ is locally complete intersection.  Our second main
result contains a description of all noetherian algebra retracts with
finitely generated homology algebra.

Let $\Max S$ denote the set of maximal ideals of $S$.

\begin{second}
Let $S\xra{\psi}R\xra{\vf}S$ be an algebra retract with noetherian
ring $R$, and set $\Max' S=\{\fn\in\Max S\var\chr(S/\fn)>0\}$.
The following conditions are equivalent.
\begin{enumerate}[{\quad\rm(i)}]
\item
The $S$-algebra $\Tor{\bu}RSS$ is finitely generated.
\item
For  every $S$-algebra $T$ there exists an isomorphism of graded
$T$-algebras
\[
\Tor{\bu}RS{T}\cong
\big({\ts\bigwedge}_{S}D_1\otimes_S\sym_{S}D_2\big)\otimes_ST
\]
where $D_1$ and $D_2$ are projective $S$-modules concentrated in degrees
$1$ and $2$, respectively, and $(D_2)_{\fn'}=0$ for all $\fn'\in\Max' S$.
\item
The $S$-modules $\aq 1SRS$ and $\aq 2SRS$ are projective, $\aq 3SRS=0$,
and $\aq 2SRS_\fn=0$ for all $\fn\in\Max' S$.
\item
For each $\fn\in\Spec S$, the ideal $\fb$ in some (respectively, every)
commutative diagram \eqref{factoriz} is generated by a regular sequence
contained in $(\bsx)^2$, and $\fb=0$ if $\fn$ is contained in some
$\fn'\in\Max' S$.
 \end{enumerate}
  \end{second}

If $S$ is a flat algebra over some ring $A$, then
$\Tor{\bu}{S\otimes_AS}SS$ is isomorphic to the {\em Hochschild homology
algebra\/} $\hch_\bu(S|A)$ of $S$ over $A$.  Our main result in \cite{hh}
shows that if the ring $R=S\otimes_AS$ is noetherian, and $\hch_\bu(S|A)$
is finitely generated as an algebra over $S$, then $S$ is regular over
$A$.  On the other hand, by the Hochschild-Kostant-Rosenberg Theorem
\cite{HKR}, as generalized by Andr\'e \cite{sym}, if $S$ is regular over
$A$ then $\hch_\bu(S|A)\cong{\ts\bigwedge}_{S}D_1$.  Thus, in the context
of Hochschild homology the module $D_2$ in Theorem II is trivial. It is
also trivial for algebra retracts where all the residue fields of $S$
have positive characteristic.  However, $\BQ\to\BQ[x]/(x^2)\to\BQ$
has finitely generated Tor algebra with $D_2\ne0$.

We proceed with an overview of the contents of the article.  Although its
main topic is the simplicially defined Andr\'e-Quillen homology theory,
many arguments are carried out in the context of DG (= differential
graded) homological algebra.

Section \ref{Differential graded algebras} contains basic definitions
and results on DG algebras.

In Section \ref{Factorizations of local homomorphisms} we recall the
construction and first properties of non-negative integers $\ve_n(\vf)$,
attached in \cite{lci} to every local homomorphism $\vf$.  These {\em
deviations\/}, whose vanishing characterizes regularity and c.i.\
properties of $\vf$, are linked to certain Andr\'e-Quillen homology
modules, but are easier to compute.  Section \ref{Indecomposables}
contains a general theorem on morphisms of minimal models of local rings.
Its proof is long and difficult.  Its applications go beyond the present
discussion.

The next two sections are at the heart of our investigation.

In Section \ref{Almost small local homomorphisms} we define a class of
local homomorphisms, that we call {\em almost small\/}.  It contains
the small homomorphisms introduced in \cite{sma}, and its larger size
offers technical advantages that are essential to our study.  We provide
various characterizations of almost small homomorphisms and give examples.
The key result established in this section is a structure theorem for
surjective almost small homomorphisms of complete rings in terms of
morphisms of DG algebras.

The proof of Theorem I depends on another new concept---that of {\em weak
category\/} of a local homomorphism.  It is defined in Section \ref{Weak
category of a local homomorphism}, where arguments from \cite{lci} are
adapted in order to obtain information on the positivity and growth of
deviations of homomorphisms with finite weak category.  To apply these
results to almost small homomorphisms we prove that they have finite
weak category; the proof involves most of the material developed up to
that point.

In Section \ref{Andre-Quillen homology} we return to Andr\'e-Quillen
homology, focusing on local homomorphisms of local rings.  We show that
vanishing of homology with coefficients in the residue field characterizes
complete intersection homomorphisms among the homomorphisms having finite
weak category.  This leads to local versions of Theorems I and II above.
The theorems themselves are proved in the final Section
\ref{Andre-Quillen dimension}.

The main results of this paper were announced in \cite{tata}, cf.\ also
Remark \ref{errata}.  That article provides historical background,
a more leisurely discussion of applications of Andr\'e-Quillen homology
to the structure of commutative algebras, and new proofs of some earlier
results on the subject.  Recently, J.\ Turner \cite{Tu} has started a
study of nilpotency in the homotopy of simplicial commutative algebras
over a field of characteristic $2$, with a view towards applications to
Quillen's Conjecture.

\section{Differential graded algebras}
\label{Differential graded algebras}

We use the theory of Eilenberg-Moore derived functors as described in
\cite[\S 1,\S 2]{AH1}.  We recall a minimum of material, referring for
details to {\em loc.\ cit\/}.

\begin{chunk}
\label{DG algebras}
Every graded object is concentrated in non-negative degrees, the
differential of every complex has degree $-1$, and each DG algebra $C$
is {\em graded commutative\/}:
\[
c'c'' = \sgn{c'}{c''} c''c' \text{ for all }c',c''\in C
\quad\text{and}\quad c^2 = 0
\text{ for all }
c\in C \text{ with $|c|$ odd}
\]
where $|c|$ denotes the degree of $c$.  The graded algebra underlying
$C$ is denoted $C\nat$.

We set $C^{[2]}=C_0+\dd(C_1)C_{\ges1}+(C_{\ges1})^2$ and $\ind C=
C/C^{[2]}$. This is a complex of $\HH0C$-modules and every morphism of
DG algebras $\gamma\col C\to D$ induces a morphism of complexes of
$\HH0D$-modules $\ind{\gamma}\col\ind{C}\otimes_{\HH0C}\HH0D\to\ind{D}$.
 \end{chunk}

\begin{chunk}
\label{quisms}
A morphism $\gamma\col C\to C'$ of DG algebras is a {\em
quasiisomorphism\/} if it induces an isomorphism in homology; this
is often signaled by the appearance of the symbol $\simeq$ next to
its arrow.  Let $C\to E$ be a morphism of $DG$ algebras, such that the
$C\nat$-module $E\nat$ is flat.  If $\gamma$ is a quasiisomorphism, then
so is $\gamma\otimes_C E \col E\to C'\otimes_CE$.  If $\epsilon\col E\to
E'$ is a quasiisomorphism and the graded $C\nat$-module $E'\nat$ is flat
as well, then $C'\otimes_C\epsilon\col C'\otimes_CE\to C'\otimes_CE'$
is a quasiisomorphism.
 \end{chunk}

\begin{chunk}
A {\em semifree extension\/} of $C$ is a DG algebra $C[X]$ such that
$C[X]\nat$ is isomorphic to the tensor product over $\BZ$ of $C\nat$
with the {\em symmetric algebra\/} of a free $\BZ$-module with basis
$\bigsqcup_{i\ges 0}X_{2i}$ and the {\em exterior algebra\/} of a free
$\BZ$-module with basis $\bigsqcup_{i\ges 0}X_{2i+1}$; the differential
of $C[X]$ extends that of $C$.

A {\em semifree \vG-extension\/} of $C$ is a DG algebra $C\la X'\ra$
such that $C\la X'\ra\nat$ is isomorphic to the tensor product over
$\BZ$ of $C\nat$ with the {\em symmetric algebra\/} of a free
$\BZ$-module with basis $X'_{0}$, the {\em exterior algebra\/} of a
free $\BZ$-module with basis $\bigsqcup_{i\ges 0}X'_{2i+1}$ and the
{\em divided powers algebra\/} of a free $\BZ$-module with basis
$\bigsqcup_{i\ges 1}X'_{2i}$; the differential of $C\la X'\ra$ extends
that of $C$, and for every $x'\in X'_{2i}$ with $i\ge1$ the $j$th
divided power $\dvip j{x'}$ satisfies $\dd(\dvip
j{x'})=\dd(x')\dvip{j-1}{x'}$ for all $j\ge1$.
 \end{chunk}

\begin{chunk}
\label{semifree}
Any morphism of DG algebras $C\to D$ factors as the canonical injection
$C\hra C[X]$ followed by a surjective quasiisomorphism $C[X]\tra D$.
If $\rho\col F\to D'$ is a surjective quasiisomorphism, then for each
commutative diagram
\[
\xymatrixrowsep{2pc}
\xymatrixcolsep{2.5pc}
\xymatrix{
\ C\
\ar@{^{(}->}[r]
\ar@{->}[d]^{\gamma}
& C[X]
\ar@{->}[r]
\ar@{-->}[d]_{\ul\delta}
& D
\ar@{->}[d]_{\delta}
\\
C'
\ar@{->}[r]
& F
\ar@{->>}[r]^\simeq_{\rho}
& D'
}
\]
of morphisms of DG algebras displayed by solid arrows there exists a
unique up to $C$-linear homotopy morphism $\ul\delta$ preserving
commutativity.
 \end{chunk}

\begin{chunk}
\label{DG Tor}
The diagrams $D\gets C\to E$ of DG algebras are the objects of a
category, whose morphisms are commutative diagrams of DG algebras
\[
\xymatrixrowsep{2pc}
\xymatrixcolsep{2.5pc}
\xymatrix{
D
\ar@{<-}[r]
\ar@{->}[d]_{\delta}
& C
\ar@{->}[r]
\ar@{->}[d]_{\gamma}
& E
\ar@{->}[d]_{\epsilon}
\\
D'
\ar@{<-}[r]
& C'
\ar@{->}[r]
& E'
}
\]
In view of \ref{semifree}, $\Tor{\bu}CDE=\hh{C[X]\otimes_CE}$ and
$\Tor{\bu}{\gamma}{\delta}{\epsilon}=\hh{\ul\delta\otimes_\gamma
\epsilon}$ define a functor from this category to that of
graded algebras.  A fundamental property of this functor is: If
$\gamma$, $\delta$, $\epsilon$ above are quasiisomorphisms, then
$\Tor{\bu}{\gamma}{\delta}{\epsilon}$ is bijective.  By \ref{quisms},
each factorization $C\to F\xra{\simeq}D$ with $F\nat$ flat over $C\nat$
yields a unique isomorphism $\Tor{\bu}CDE\to \hh{F\otimes_CE}$ of graded
algebras.
 \end{chunk}

\begin{chunk}
\label{Gamma algebra}
A {\em DG \vG-algebra\/} is a DG algebra $K$ in which a sequence
$\{\dvip jx\in K_{jn}\}_{j\ges 0}$ of {\em divided powers\/} is defined
for each $x\in K_n$ with $n$ even positive, and satisfies a list of
standard identities; it can be found in full, say, in \cite[(1.7.1),
(1.8.1)]{GL}.  A {\em morphism of DG \vG-algebras\/} $\varkappa\col
K\to L$ is a morphism of DG algebras such that $\varkappa(\dvip jx)=
(\varkappa(x))^{(j)}$ for all $x\in K$ with $|x|$ even positive and all
$j\in\BN$.

Let $K^{(2)}$ denote the $K_0$-submodule of $K$ generated by $K^{[2]}$
and all $\dvip jx$, where $|x|$ is even positive and $j\geq 2$.  Set $\gin
K=K/K^{(2)}$.  This is a complex of $\HH0K$-modules.  Every morphism of
\G-algebras $\varkappa\col K\to L$ induces a morphism $\gin{\varkappa}
\col \gin{K}\otimes_{\HH0K}\HH 0L\to\gin{L}$ of complexes of $\HH
0L$-modules.
 \end{chunk}

\begin{chunk}
\label{Gamma semifree}
If $K$ is a DG \G-algebra, then $K\la X'\ra$ has a unique structure of
DG \G-algebra extending that of $K$ and preserving the divided powers
of the variables $x'\in X'_{2i}$ with $i>0$.  Every morphism of DG
\vG-algebras $\varkappa\col K\to L$  can be factored as $K\hra K\la
X'\ra\tra L$ with second map a surjective quasiisomorphism of DG
\G-algebras. If $\zeta\col M\to L'$ is a surjective quasiisomorphism,
then for each commutative diagram
\[
\xymatrixrowsep{2pc}
\xymatrixcolsep{2.5pc}
\xymatrix{
\ K\
\ar@{^{(}->}[r]
\ar@{->}[d]_{\varkappa}
& K\la X'\ra
\ar@{->}[r]
\ar@{-->}[d]_{\ul\lambda}
& L
\ar@{->>}[d]_{\lambda}
\\
K'
\ar@{->}[r]
& M
\ar@{->>}[r]^\simeq_{\zeta}
& L'
}
\]
of morphisms of DG algebras displayed by solid arrows there exists
a unique up to $K$-linear homotopy morphism of DG \G-algebras
$\ul\lambda$ making both squares commute.
 \end{chunk}

\begin{chunk}
\label{DG Gamma Tor}
Divided powers of a cycle are cycles, but divided powers of a boundary
need not be boundaries. If they are, then the DG \G-algebra $K$ is
called {\em admissible\/}, and $\hh K$ inherits from $K$ a structure of
\G-algebra. This notion of admissibility is less restrictive than the
one adopted in \cite{AH1}, and lacks some of the desirable properties
the latter posesses, but it suffices for the needs of this paper.

Let $K\tra k$ be a surjective morphism of DG algebras, where $k$ is a
field concentrated in degree $0$, and let $k\hra l$ be a field
extension.  If $K\hra K\la X'\ra\tra k$ is a factorization as in
\ref{Gamma semifree}, then the unique DG \G-algebra structure on $l\la
X'\ra= K\la X'\ra\otimes_Kl$ is admissible, cf.\ \cite[(2.6)]{AH1} or
\cite[(3.4)]{tata}.  Thus, Tor defines a functor from the category of
diagrams $K\tra k\hra l$, with the obvious morphisms, to the category
of \G-algebras and their morphisms.
 \end{chunk}

\section{Factorizations of local homomorphisms}
\label{Factorizations of local homomorphisms}

Let $\vf \col (R,\fm,k)\to (S,\fn,l)$ be a homomorphism of local rings,
which is {\em local\/} in the sense that $\vf(\fm)\subseteq\fn$.  A
{\em regular factorization\/} of $\vf$ is a commutative diagram
\[
\xymatrixrowsep{2pc}
\xymatrixcolsep{3pc}
\xymatrix{
&R'
\ar@{->}[dr]
\\
R
\ar@{->}[rr]^{\vf}
\ar@{->}[ur]
&& S
}
\]
of local homomorphisms such that the $R$-module $R'$ is flat, the ring
$R'/\fm R'$ is regular, and the map $R'\to S$ is surjective.

Regular factorizations are often easily found, for instance, when $\vf$
is essentially of finite type (in particular, surjective), or when
$\vf$ is the canonical embedding of $R$ in its completion with respect
to the maximal ideal.  In this paper they are mostly used through the
following construction of Avramov, Foxby, and B.\ Herzog \cite{afh}.

\begin{chunk}
\label{cohen}
If $\grave\vf\col R\to \wh S$ is the composition of $\vf$ with the
canonical inclusion $S\to\wh S$, then by \cite[(1.1)]{afh}, $\grave\vf$
has a regular factorization $R\xra{\dot\vf}R'\xra{\vf'}\wh S$ with a
complete local ring $R'$; it is called a {\em Cohen factorization\/} of
$\grave\vf$.  By \cite[(1.5)]{afh}, it can be chosen to satisfy the
additional condition $\edim R'/\fm R'=\edim S/\fm S$; we say that such
a Cohen factorization is {\em reduced\/} (it is called minimal in
\cite{afh}).  Clearly, any regular factorization $\vf=\pi\circ\iota$
gives rise to a Cohen factorization $\grave\vf=\wh\pi\circ\grave\iota$.

Cohen factorizations need not be isomorphic.  However, if $R\xra{\ddot\vf}
R''\xra{\vf''}\wh S$ also is a Cohen factorization of $\grave\vf$, then by
\cite[(1.2)]{afh} there exists a commutative diagram
\[
\xymatrixrowsep{2pc}
\xymatrixcolsep{3pc}
\xymatrix{
& R'
\ar@{->}[dr]^{\varphi'}
\\
R
\ar@{->}[ur]^{\dot\varphi}
\ar@{->}[r]
\ar@{->}[dr]_{\ddot\varphi}
& R'''
\ar@{->}[u]
\ar@{->}[r]
\ar@{->}[d]
& {\widehat S}
\\
& R''
\ar@{->}[ur]_{\varphi''}
}
\]
of local homomorphisms, where the horizontal row is a Cohen factorization,
and the vertical maps are surjections with kernels generated by regular
sequences whose images in $R'''/\fm R'''$ can be completed to regular
systems of parameters.
 \end{chunk}

\begin{chunk}
\label{decomposable}
Let $(A,\fp,k)$ be a local ring. We say that a semifree extension $A[X]$
has {\em decomposable differential\/} if $X=X_{\ges1}$ and
\[
\dd(X)\subseteq\fp A[X]+(X)^2A[X]\,.
\]
When this condition holds, for each $n\ge 1$ there are equalities
\[
\HH n{A[X]/(\fp,X_{\sles n})}=\ZZ_n(A[X]/(\fp,X_{\sles n}))=kX_n\,.
\]
 \end{chunk}

\begin{chunk}
\label{models}
Let $\vf\col (R,\fm,k)\to(S,\fn,l)$ be a local homomorphism.

A {\em minimal model\/} of $\vf$ is a diagram $R\xra{\wt\iota}R'[U]
\xra{\wt\vf} S$ where the differential of $R'[U]$ is decomposable,
$\wt\vf$ is a quasiisomorphism, and $\vf=\wt\vf_0\circ\wt\iota_0$
is a regular factorization.  If $\vf$ has a regular factorization
(in particular, if $S=\wh S$), then $\vf$ has a minimal model:  The
DG algebra $R'[U]$ is obtained by successively adjoining to $R'$ sets
of variables $U_n$ of degree $n\ge1$, so that $\dd(U_1)$ minimally
generates $\Ker(\pi)$ and $\dd(U_n)$ is a minimal set of generators for
$\HH{n-1}{R[U_{<n}]}$, cf.\ \cite[(2.1.10)]{res}.
 \end{chunk}

The next proposition elaborates on \cite[(3.1)]{lci}.

\begin{proposition}
\label{independence}
Let $\vf\col(R,\fm,k)\to(S,\fn,l)$ be a local homomorphism and let
$R\to R'[U']\to \wh S$ and $R\to R''[U'']\to \wh S$ be minimal models
of $\grave\vf$.

For each integer $n\ge2$ there are equalities
\begin{align*}
\card(U_1') - \edim(R'/\fm R')&= \card(U_1'') - \edim(R''/\fm R'')\,,\\
\card(U_n')&= \card(U_n'')\,,
\end{align*}
and there exist isomorphisms of DG algebras over the field $l$
\[
l[U'_{\ges n}]=R'[U']/(\fm',U'_{< n})\cong
R''[U'']/(\fm'',U''_{< n})= l[U''_{\ges n}]\,.
\]
 \end{proposition}

\begin{proof}
By \ref{cohen} we may assume there is a surjection $R''\to R'$ with
kernel generated by a regular sequence $\bsx$ that extends to a minimal
generating set of the maximal ideal $\fm''$ of $R''$.  Changing $U_1''$
if need be, we may assume that $U''=V\sqcup U$ with $\dd(V)=\bsx$.
The canonical map $R''[V]\to R'$ is a quasiisomorphism, $R''[V]$ is a
DG subalgebra of $R''[U'']$ and the $R''[V]\nat$-module $R''[U'']\nat$
is free, so the induced map
\[
R''[U''] \to R''[U'']/(V,\dd(V)) =R'[U]
\]
is a quasiisomorphism, cf.\ \ref{quisms}.  Thus, $\hh{R'[U]} \cong \wh
S$, and the differential of $R'[U]$ is decomposable because it is
induced by that of $R''[U'']$.  By \cite[(7.2.3)]{res} there exists an
isomorphism $R'[U']\cong R'[U]$ of DG algebras over $R'$, so we get
\[
R'[U']/(\fm',U'_{<n})\cong R'[U]/(\fm',U_{<n})
\]
for all $n\ge 1$.  The algebra on the right is equal to
$R''[U'']/(\fm'',U''_{<n})$ for $n\ge2$, so we have proved the last
assertion.  In view of \ref{decomposable}, it implies
\[
lU'_n=\HH n{R'[U']/(\fm',U'_{<n})}\cong \HH n{R'[U]/(\fm',U_{<n})}=lU_n\,.
\]
Thus, we obtain numerical equalities
\begin{align*}
\card(U'_n)&=\card(U_n)=\card(U''_n)\qquad\text{for}\qquad n\ge2\,;\\
\card(U_1')&=\card(U_1)=\card(U_1'')-\card(V)\\
&=\card(U_1'')-\big(\edim(R''/\fm R'')-\edim(R'/\fm R')\big)\,.
\end{align*}

All the assertions of the proposition have now been established.
\end{proof}

\begin{definition}
\label{deviations}
Let $\vf\col(R,\fm,k)\to S$ be a local homomorphism, and let $R\to
R'[U]\to \wh S$ be a minimal model of $\grave\vf$.  The $n$th {\em
deviation\/} of $\vf$ is the number
\[
\ve_n(\vf) =
\begin{cases}
\card(U_1) - \edim(R'/\fm R')+\edim(S/\fm S)
& \text{ for } n=2\,;\\
\card(U_{n-1})
& \text{ for } n\ge 3\,.
  \end{cases}
\]
By Proposition \ref{independence}, these are invariants of $\vf$.
Deviations were defined in \cite[\S 3]{lci} with a typo in the
expression for $\ve_2(\vf)$, which is corrected above.

Note that $\ve_n(\vf)\ge0$ for all $n$: this is clear for $n\ge3$;
for $n=2$, use the equalities
\begin{align*}
\card(U_1)
&=\rank_l\left(\frac{\Ker(\vf')}{\fm'\Ker(\vf')}\right)\,;\\
\edim(R'/\fm R')-\edim(S/\fm S)
&=\rank_l\left(\frac{\Ker(\vf')}{\Ker(\vf')\cap(\fm'{}^2+\fm R')}\right)\,.
\end{align*}
 \end{definition}

Vanishing of deviations is linked to the structure of $\vf$.  We
reproduce \cite[(3.2)]{lci}:

\begin{proposition}
\label{regular map}
If $\vf\col(R,\fm,k)\to S$ is a local homomorphism, then the following
conditions are equivalent.
\begin{enumerate}[\rm\quad(i)]
\item
$\vf$ is flat and $S/\fm S$ is regular.
\item
$\ve_n(\vf)=0$ for all $n\ge 2$\,.
\item
$\ve_2(\vf)=0$\,.
\end{enumerate}
\end{proposition}

\begin{proof}
(i) $\implies$ (ii)
The diagram $R\to\wh S=\wh S$ is a Cohen factorization of $\grave\vf$,
so $\grave\vf$ has a minimal model with $U=\varnothing$.

(iii) $\implies$ (i)
Choose a reduced Cohen factorization. By definition, $\ve_2(\vf)=0$
entails $U_1=\varnothing$, so $\wh S=\HH0{R'[U]}=R'$, hence $\wh S$ is
flat over $R$ and $\wh S/\fm \wh S$ is regular; these properties
descend to $S$ and $S/\fm S$.
 \end{proof}

The following notion is basic for the rest of the paper.

\begin{definition}
\label{def: ci}
A local homomorphism $\vf\col R\to (S,\fn,l)$ is {\em complete
intersection\/} (or {\em c.i.\/}) at $\fn$, if in some Cohen
factorization $R\to R'\xra{\vf'}\wh S$ of $\grave\vf$ the ideal
$\Ker(\vf')$ is generated by an $R'$-regular sequence.
 \end{definition}

Other definitions of c.i.\ homomorphisms require additional hypotheses
on $\vf$; when they hold, the general concept specializes properly,
cf.\ \cite[(5.2), (5.3)]{lci}. The next proposition amplifies
\cite[(3.3)]{lci}; it shows, in particular, that the c.i.\ property is
detected by every Cohen factorization.

\begin{proposition}
\label{ci map}
If $\vf\col R\to(S,\fn,l)$ is a local homomorphism, then the following
conditions are equivalent.
\begin{enumerate}[\rm\quad(i)]
\item
$\vf$ is complete intersection at $\fn$\,.
\item
$\ve_n(\vf)=0$ for all $n\ge 3$\,.
\item
$\ve_3(\vf)=0$\,.
\end{enumerate}
\end{proposition}

\begin{proof}
In any minimal model $R\to R'[U]\to\wh S$ of $\grave\vf$ the DG
algebra $R'[U_1]$ is the Koszul complex on a minimal set of generators
of $\Ker(\vf')$.  If (i) holds, then $U=U_1$, so (i) implies (ii).
If (iii) holds, then $\HH 1{R'[U_1]}=0$, so the ideal $\Ker(\vf')$
is generated by a regular sequence.
 \end{proof}

\section{Indecomposables}
\label{Indecomposables}

In this section we analyze the divided powers in Tor.

\begin{chunk}
If $(R,\fm,k)$ is a local ring, then $\Tor{\bu}{R}kk$ is a \G-algebra,
cf.\ \ref{DG Gamma Tor}.

Using the functor $\gin{-}$ of \G-indecomposables defined in \ref{Gamma
algebra}, we set
\[
\pies R=\gin{\Tor{\bu}{R}kk}
\]
If $\ov\vf\col k\rat l$ is a field extension, then the canonical
isomorphism
\[
\Tor{\bu}Rkk\otimes_kl\cong\Tor{\bu}Rkl
\]
is one of \G-algebras, and so induces an isomorphism of graded
$l$-vector spaces
\[
\pies R\otimes_kl\cong\gin{\Tor{\bu}{R}kl}
\]
that we use as identification.  Thus, every local homomorphism $\vf \col
R\to(S,\fn,l)$ defines an $l$-linear homomorphism of graded vector spaces
\[
\pies\vf\col\pies R\otimes_kl \xra{\gin{\Tor{\bu}{\vf}{\ov\vf}l}}\pies S\,.
\]
 \end{chunk}

\begin{example}
\label{acyclic closure}
Let $(R,\fm,k)$ be a local ring.  An {\em acyclic closure\/} of $k$ is
a factorization $R\to R\la X'\ra \to k$ of the epimorphism $R\to k$, as
in \ref{Gamma semifree}, constructed so that $\dd(X'_1)$ minimally
generates $\fm$ and $\dd(X'_n)$ minimally generates $\HH{n-1}{R\la
X'_{<n}\ra}$ for each $n\geq 2$, cf.\ \cite[(6.3)]{res}.  By an
important theorem of Gulliksen \cite{Gul} and Schoeller \cite{Sch}, in
this case $\dd(R\la X'\ra)\subseteq \fm R\la X'\ra$, cf.\ also
\cite[(6.3.4)]{res}.  This yields isomorphisms
\[
\pie nR\cong kX'_n \qquad\text{for all}\qquad n\in\BZ\,.
\]

The $n$th {\em deviation\/} of $R$ is the number $\ve_n(R)=\card{X'_n}$.
They measure the singularity of $R$: $\ve_n(R)=0$ for all $n\ge2$ if
and only if $\ve_2(R)=0$, if and only if $R$ is regular; $\ve_n(R)=0$
for all $n\ge3$ if and only if $\ve_3(R)=0$, if and only if $R$ is
c.i., cf.\ \cite[Ch.\ III]{GL}, \cite[\S 7]{res}. These results can be
derived from Propositions \ref{regular map} and \ref{ci map}, since by
\cite[(7.2.5)]{res} deviations of rings and of homomorphisms are linked
as follows:
 \end{example}

\begin{chunk}
\label{comparison}
If $\vf\col A\to R$ is a surjective local homomorphism with $A$
regular, then
\[
\ve_n(\vf)=\ve_n(R) \quad\text{for all}\quad n\ge2\,.
\]
 \end{chunk}

The next result is a functorial enhancement of the numerical equality
above.

\begin{theorem}
\label{indecomposables}
Consider a commutative diagram of morphisms of DG algebras
\begin{equation}
\label{original}
\begin{aligned}
\xymatrixrowsep{2pc}
\xymatrixcolsep{2pc}
\xymatrix{
A_{\mathstrut}
\ar@{->}[r]^{\beta}
\ar@/_2pc/[dd]_{\rho}
\ar@{^{(}->}[d]
&
B_{\mathstrut}
\ar@/^2pc/[dd]^{\sigma}
\ar@{_{(}->}[d]
\\
A[X]
\ar@{->}[r]^{\phi}
\ar@{->>}[d]_{\wt\rho}^{\simeq}
& B[Y]
\ar@{->>}[d]^{\wt\sigma}_{\simeq}
\\
R
\ar@{->}[r]^{\vf}
& S
}
\end{aligned}
\end{equation}
where $(R,\fm,k)$ and $(S,\fn,l)$ are local rings, $(A,\fp,k)$
and $(B,\fq,l)$ are regular local rings, the homomorphisms $\vf$ and
$\beta$ are local, the homomorphisms $\rho$ and $\sigma$ are surjective,
$\Ker(\rho)\subseteq\fp^2$ and $\Ker(\sigma)\subseteq \fq^2$, and the
triangles are minimal models.

For each $n\ge2$ there exists a commutative diagram of homomorphisms
\[
\xymatrixrowsep{2pc}
\xymatrixcolsep{5pc}
\xymatrix{
\pie{n}{R}\otimes_k l
\ar@{->}[d]_{\cong}
\ar@{->}[r]^-{\pie{n}{\vf}}
& \pie{n}{S}
\ar@{->}[d]^{\cong}
\\
\Ind{n-1}{l[X]}
\ar@{->}[r]^-{\Ind{n-1}{\phi\otimes_\beta l}}
& \Ind{n-1}{l[Y]}
}
\]
of $l$-vector spaces, where the vertical arrows are isomorphisms.
 \end{theorem}

The theorem shows that $\pie{n}{\vf}$ and $\Ind{n-1}{\phi\otimes_\beta l}$
determine each other.  These are very different maps: the first is induced
by a morphism of DG \G-algebras, while divided powers have no role in
the construction of the second.  This accounts for the intricacies of the
proof.  In it, and later in the paper, it is convenient to suppress the
effect of $\Tor1Rkk$ on $\Tor{\bu}Rkk$.  We do that in a systematic way.

\begin{chunk}
\label{rtor} The {\em reduced torsion algebra\/} of a local ring
$(R,\fm,k)$ is the $k$-algebra
\[
\rTor{\bu}Rkk = \frac{\Tor{\bu}Rkk}{\Tor{\bu}Rkk\cdot\Tor1Rkk}\,.
\]

Since $\Tor{\bu}Rkk$ is a \G-algebra and the ideal
$J=\Tor{\bu}Rkk\cdot\Tor1Rkk$ is generated by elements of degree $1$,
basic properties of divided powers imply that each element of even
degree $a\in J$ satisfies $\dvip ia\in J$ for all $i\ge1$.  It follows
that $\rTor{\bu}Rkk$ admits a unique \G-structure for which the
canonical surjection $\Tor{\bu}Rkk\to \rTor{\bu}Rkk$ becomes a morphism
of \G-algebras, hence
\[
\pie{\ges2}R=\gin{\rTor{\bu}Rkk}\,.
\]
If $\vf\col R\to (S,\fn,l)$ is a local homomorphism, then
$\Tor{\bu}{\vf}{\ov\vf}l$ induces a morphism
\[
\rTor{\bu}{\vf}{\ov\vf}l\col \rTor{\bu}Rkl\to\rTor{\bu}Sll
\]
\G-algebras, so for $n\ge2$ we get commutative diagrams of $l$-linear
homomorphisms
\[
\xymatrixrowsep{2pc}
\xymatrixcolsep{6pc}
\xymatrix{
\pie{n}{R}\otimes_k l
\ar@{->}[d]_{\cong}
\ar@{->}[r]^-{\pie{n}{\vf}}
& \pie{n}{S}
\ar@{->}[d]^{\cong}
\\
\Gin{n}{\rTor{\bu}Rkl}
\ar@{->}[r]^-{\Gin n{\rTor{\bu}{\vf}{\ov\vf}l}}
&
\Gin{n}{\rTor{\bu}Sll}
}
\]
\end{chunk}

The {\em proof\/} of Theorem \ref{indecomposables} takes up the rest of
the section.  Only its statement is used later, so the reader may skip
to the next section without loss of continuity.

We start by forming a diagram of morphisms of DG \G-algebras
\begin{equation}
\label{koszuls}
\begin{aligned}
\xymatrixrowsep{2pc}
\xymatrixcolsep{2pc}
\xymatrix{
A_{\mathstrut}
\ar@{->}[r]^{\beta}
\ar@{^{(}->}[d]
&
B_{\mathstrut}
\ar@{_{(}->}[d]
\\
A\la X'_1\ra
\ar@{->}[r]^{\varkappa}
\ar@{->>}[d]_{\epsilon}^{\simeq}
& B\la Y'_1\ra
\ar@{->>}[d]^{\eta}_{\simeq}
\\
k
\ar@{->}[r]^{\ov\vf}
& l
}
\end{aligned}
\end{equation}
in the following order. First we form the vertical sides by choosing them
to be acyclic closures of the respective residue fields.  Next we note
that since both $A$ and $B$ are regular local rings, the DG algebras
$A\la X'_1\ra$ and $B\la Y'_1\ra$ are Koszul complexes on minimal
sets of generators of $\fp$ and $\fq$, respectively.  Finally, we use
\ref{Gamma semifree} to choose a morphism $\varkappa$ that preserves
the commutativity of the rectangle.

Base change from Diagram \eqref{koszuls} yields the central
rectangles in the
diagram
\begin{equation}
\label{closures}
\begin{aligned}
\xymatrixrowsep{2pc}
\xymatrixcolsep{2pc}
\xymatrix{
R\la X'\ra
\ar@/^2pc/[rrr]^{\ul\vf}
\ar@{=}[d]
& R_{\mathstrut}
\ar@{->}[r]^{\vf}
\ar@{_{(}->}[l]
\ar@{^{(}->}[d]_{\tau}
&
S_{\mathstrut}
\ar@{^{(}->}[r]
\ar@{_{(}->}[d]^{\theta}
& S\la Y'\ra
\ar@{=}[d]
\\
R\la X'\ra
\ar@{->>}[rd]_{\simeq}
& R\la X'_1\ra
\ar@{_{(}->}[l]
\ar@{->}[r]^{\vf\otimes_\beta\varkappa}
\ar@{->>}[d]
& S\la Y'_1\ra
\ar@{^{(}->}[r]
\ar@{->>}[d]
& S\la Y'\ra
\ar@{->>}[ld]^{\simeq}
\\
& k
\ar@{->}[r]^{\ov\vf}
& l
}
\end{aligned}
\end{equation}
of morphisms of DG \G-algebras.  The rest is constructed as follows.
In view of the hypotheses on $\rho$ and $\sigma$, minimal sets of
generators of $\fp$ and $\fq$ map to minimal sets of generators of $\fm$
and $\fn$, respectively.  By Example \ref{acyclic closure} the DG algebras
$R\la X'_1\ra= R\otimes_AA\la X'_1\ra$ and $S\la Y'_1\ra=S\otimes_BB\la
Y'_1\ra$ can be extended to acyclic closures $R\hra R\la X'\ra\tra k$
and $S\hra S\la Y'\ra\tra l$.  Finally, the morphism $\ul\vf$ is chosen
so as to preserve the commutativity of the diagram: this is possible by
\ref{Gamma semifree}.

\begin{lemma}
\label{reduced}
Diagram {\em\eqref{closures}} induces a commutative diagram
\[
\begin{aligned}
\xymatrixrowsep{2pc}
\xymatrixcolsep{1.3pc}
\xymatrix{
l\la X'_{\ges2}\ra
\ar@{->}[r]^-{\cong}
&\rTor{\bu}Rkl
\ar@{->}[rrrr]^{\rTor{\bu}{\vf}{\ov\vf}l}
\ar@{->}[d]_{\cong}
&&&& \rTor{\bu}Sll
\ar@{->}[d]^-{\cong}
& l\la Y'_{\ges2}\ra
\ar@{->}[l]_{\cong}
\\
& \Tor{\bu}{R\la X'_1\ra}kl
\ar@{->}[rrrr]^{\Tor{\bu}{\vf\otimes_\beta\varkappa}{\ov\vf}l}
&&&& \Tor{\bu}{S\la
Y'_1\ra}ll
}
\end{aligned}
\]
of homomorphisms of \vG-algebras.
 \end{lemma}

\begin{proof}
By construction, $R\la X'\ra$ and $S\la Y'\ra$ are acyclic closures.
In view of \ref{acyclic closure}, this means that there are inclusions
$\dd(R\la X'\ra)\subseteq\fm R\la X'\ra$ and $\dd(S\la Y'\ra)\subseteq
\fn S'\la Y'\ra$.  These inclusions provide the equalities in the
commutative diagram
\[
\xymatrixrowsep{2pc}
\xymatrixcolsep{1pc}
\xymatrix{
l\la X'\ra
\ar@{=}[r]
\ar@{->}[d]
\ar@/^2pc/[rrrrr]^{\ul\vf\otimes_{\vf}l}
& \Tor{\bu}Rkl
\ar@{->}[rrr]^{\Tor{\bu}{\vf}{\ov\vf}l}
\ar@{->}[d]_{\Tor{\bu}\tau kl}
&&& \Tor{\bu}Sll
\ar@{=}[r]
\ar@{->}[d]^{\Tor{\bu}\theta ll}
& l\la Y'\ra
\ar@{->}[d]
\\
l\la X'_{\ges2}\ra
\ar@{=}[r]
& \Tor{\bu}{R\la X'_1\ra}kl
\ar@{->}[rrr]^{\Tor{\bu}{\vf\otimes_\beta\varkappa}{\ov\vf} l}
&&& \Tor{\bu}{S\la
Y'_1\ra}ll
\ar@{=}[r]
& l\la Y'_{\ges2}\ra
}
\]
induced by Diagram \eqref{closures}.  By \ref{DG Gamma Tor}, all the
maps are morphisms of \G-algebras.

The inclusions noted above also show that the external vertical
maps are the canonical surjections of graded algebras, whose
kernels are the ideals generated by $X'_1$ and $Y'_1$ respectively.
It follows that $\Ker(\Tor{\bu}\tau kl)$ is generated by $\Tor1Rkl$,
and $\Ker(\Tor{\bu}\theta ll)$ is generated by $\Tor1Sll$.  In view of
the definition of the reduced Tor functor in \ref{rtor}, the diagram
above induces the desired diagram.
 \end{proof}

We refine Diagram {\em\eqref{original}} to a commutative diagram of
morphisms of DG algebras
\begin{equation}
\label{refined}
\begin{aligned}
\xymatrixrowsep{2pc}
\xymatrixcolsep{3pc}
\xymatrix{
A[X]
\ar@{->}[rrr]^{\phi}
\ar@{->>}[dr]^(0.65){\wt\rho}_(0.65){\simeq}
\ar@{-->}[ddr]_{\kappa}
&&& B[Y]
\ar@{->>}[dl]_(0.65){\wt\sigma}^(0.65){\simeq}
\ar@{-->}[ddl]^{\lambda}
\\
& R
\ar@{->}[r]^{\vf}
& S
\\
& A\la V\ra_{\mathstrut}
\ar@{->}[r]^{\varPhi}
\ar@{->>}[u]_{\breve\rho}^{\simeq}
&
B\la W\ra
\ar@{->>}[u]^{\breve\sigma}_{\simeq}
\\
A^{\mathstrut}
\ar@{^{(}->}[uuu]
\ar@{->}[rrr]^{\beta}
\ar@{^{(}->}[ur]
&&&
B^{\mathstrut}
\ar@{_{(}->}[uuu]
\ar@{_{(}->}[ul]
}
\end{aligned}
\end{equation}
by performing the following steps.  First we invoke \ref{Gamma semifree}
to construct factorizations $A\hra A\la V\ra\xra{\breve\rho}R$ of $\rho$
and $B\hra B\la W\ra\xra{\breve\sigma}S$ of $\sigma$.  Next we choose
by \ref{Gamma semifree} a morphism of DG \G-algebras $\varPhi$ so as to
preserve the commutativity of the already constructed part of the diagram.
Finally, we use \ref{semifree} to obtain morphisms of DG algebras $\kappa$
and $\lambda$ which preserve the commutativity of the lateral trapezoids.

It should be noted at this point that, in general, $\varPhi\kappa\ne
\lambda\phi$.  Using Diagrams \eqref{koszuls} and \eqref{refined} we
produce a diagram of morphisms of DG algebras
\begin{equation}
\label{big}
\xymatrixrowsep{3.0pc}
\xymatrixcolsep{3.8pc}
\begin{aligned}
\xymatrix{
A[X]\la X'_1\ra
\ar@{->}[rrr]^{\phi\otimes_\beta\varkappa}
\ar@{->>}[dr]^{\wt\rho\otimes_A A\la X'_1\ra}_{\simeq}
\ar@{-->}[ddr]_(.6){\kappa\otimes_A A\la X'_1\ra}^(.55){\simeq}
\ar@{->}[dddd]^{A[X]\otimes_A\epsilon}_{\simeq}
&&& B[Y]\la Y'_1\ra
\ar@{->>}[dl]_{\wt\sigma\otimes_B B\la Y'_1\ra}^{\simeq}
\ar@{-->}[ddl]^(.6){\lambda\otimes_B B\la Y'_1\ra}_(.55){\simeq}
\ar@{->}[dddd]_{B[Y]\otimes_B\eta}^{\simeq}
\\
& R\la X'_1\ra
\ar@{->}[r]^{\vf\otimes_\beta\varkappa}
& S\la Y'_1\ra
\\
& A\la V\ra\la X'_1\ra
\ar@{->>}[u]_{\breve\rho\otimes_A A\la X'_1\ra}^{\simeq}
\ar@{->}[r]_{\varPhi\otimes_\beta\varkappa}
\ar@{->>}[d]_{A\la
V\ra\otimes_A{\epsilon}}^{\simeq}
& B\la W\ra\la Y'_1\ra
\ar@{->>}[u]^{\breve\sigma\otimes_B B\la Y'_1\ra}_{\simeq}
\ar@{->>}[d]^{B\la
W\ra\otimes_B\eta}_{\simeq}
\\
& k\la V\ra
\ar@{->}[r]^{\varPhi\otimes_\beta\ov\vf}
& l\la W\ra
\\
k[X]
\ar@{-->}[ur]^{\kappa\otimes_Ak}_{\simeq}
\ar@{->}[rrr]^{\phi\otimes_\beta\ov\vf}
&&& l[Y]
\ar@{-->}[ul]_{\lambda\otimes_Bl}^{\simeq}
}
\end{aligned}
\end{equation}
where the central rectangles are formed by morphisms of DG \G-algebras,
all non-horizontal arrows are quasiisomorphisms due to \ref{quisms},
and almost all paths commute---the possible exception being the paths
around the two trapezoids with horizontal bases and hyphenated sides.

 From \ref{DG Tor} and \ref{DG Gamma Tor} we deduce the following result.

\begin{lemma}
\label{gammas}
The maps in Diagram {\em\eqref{big}} induce a commutative diagram
\begin{equation*}
\xymatrixrowsep{2.0pc}
\xymatrixcolsep{5.0pc}
\xymatrix{
\Tor{\bu}{R\la X'_1\ra}kl
\ar@{->}[d]_{\cong}
\ar@{->}[r]^{\Tor{\bu}{\vf\otimes_\beta\varkappa}{\ov\vf}l}
&
\Tor{\bu}{S\la Y'_1\ra}ll
\ar@{->}[d]^{\cong}
\\
\Tor{\bu}{k\la V\ra}kl
\ar@{->}[r]^{\Tor{\bu}{\varPhi\otimes_\beta\ov\vf}{\ov\vf}l}
&
\Tor{\bu}{l\la W\ra}ll
}
\end{equation*}
of homomorphisms of \vG-algebras, where the vertical maps are
isomorphisms. \qed
 \end{lemma}

We pause to recall some classical material on bar-constructions.

\begin{chunk}
\label{bar}
Let $C$ be a {\em connected\/} DG algebra over the field $k$, which
means that $C_0=k$ and $\dd(C_1)=0$.  The {\em bar construction\/}
$(\bark kC,\bar\dd)$ is a connected DG \G-algebra over $k$, with
multiplication (called {\em shuffle product\/}) and divided powers
constructed in \cite[Exp.\ 7, \S 1)]{Ca}; cf.\ also \cite[Ch.\ X, \S
12]{Mc}. It has a basis consisting of symbols $[c_1|c_2|\cdots|c_p]$ of
degree $p+|c_1|+\cdots +|c_{p}|$, where the $c_i$ range independently
over a basis of $C_{\ges 1}$ and $p=0,1,2\dots$.  The element
$[c_1|c_2|\cdots|c_p]$ has {\em weight\/} $p$; the weight of $x\cdot
y$ is the sum of those of $x$ and $y$; if $|x|$ is even positive, then
the weight of $\dvip ix$ is $i$ times that of $x$.  In general, the DG
\G-algebra $\bark kC$ is not admissible.

There exists a DG algebra $(\barr kC,\dd)$ such that $\barr kC \nat=
C\nat\otimes_k\bark kC\nat$ as graded algebras, $\dd$ extends the
differential of $C$, the isomorphism $\barr kC \otimes_Ck\cong\bark
kC$ is one of DG algebras, and the augmentation $\barr kC\to k$ is a
quasiisomorphism of DG algebras.  If $C$ is a DG \G-algebra, then by
\cite[Exp.\ 7, \S 5)]{Ca} so is $\barr k C$, the map $\barr kC\to\bark
kC$ is a morphism of DG \G-algebras, and $\bark kC$ is admissible.

The bar construction is natural for morphisms $\gamma\col C\to C'$ of
connected DG $k$-algebras; a morphism of DG \G-algebras $\bark
k{\gamma}\col\bark kC\to\bark k{C'}$ is given by
\begin{equation}
\label{naturality}
\bark k{\gamma}([c_1|\cdots|c_p])=[\gamma(c_1)|\cdots|\gamma(c_p)]
\end{equation}

There is a canonical isomorphism $\bark kC\otimes_kl\cong \bark
l{C\otimes_kl}$ of DG \G-algebras over $l$.  In conjunction with \ref{DG
Tor}, it induces isomorphisms of graded algebras
\[
\Tor{\bu}{C}kl\cong\hh{\barr kC\otimes_Cl}=\hh{\bark kC\otimes_kl}
                             \cong \hh{\bark l{C\otimes_kl}}
\]
which are natural with respect to morphisms of connected DG
$k$-algebras. When $C$ is a DG \G-algebra the isomorphisms above are of
\G-algebras, cf.\ \ref{DG Gamma Tor}.
 \end{chunk}

\begin{lemma}
\label{admissibles}
The DG \vG-algebras $\bark l{l[X]}$ and $\bark l{l[Y]}$ are admissible,
and the maps in Diagram {\em\eqref{big}} induce a commutative diagram
\[
\xymatrixrowsep{2pc}
\xymatrixcolsep{4pc}
\xymatrix{
\Tor{\bu}{k\la V\ra}kl
\ar@{->}[r]^{\Tor{\bu}{\varPhi\otimes_\beta\ov\vf}{\ov\vf}l}
\ar@{->}[d]_{\cong}
&
\Tor{\bu}{l\la W\ra}ll
\ar@{->}[d]^{\cong}
\\
\hark l{l[X]}
\ar@{->}[r]^{\hark l{\phi\otimes_\beta l}}
& \hark l{l[Y]}
}
\]
of morphisms of \vG-algebras over $l$, where the vertical maps are
isomorphisms.
 \end{lemma}

\begin{proof}
Diagram \eqref{big} induces a diagram of homomorphisms of graded algebras
\[
\xymatrixrowsep{3.0pc}
\xymatrixcolsep{0.4pc}
\xymatrix{
\Tor{\bu}{A[X]\la X'_1\ra}kl
\ar@{->}[rrrrr]
\ar@{->}[dr]^{\cong}
\ar@{->}[ddr]_{\cong}
\ar@{->}[dddd]^(.6){\Tor{\bu}{A[X]\otimes_A{\epsilon}}kl}_(.6){\cong}
&&&&&
\Tor{\bu}{B[Y]\la Y'_1\ra}ll
\ar@{->}[dl]_{\cong}
\ar@{->}[ddl]^{\cong}
\ar@{->}[dddd]^(.6){\cong}
\\
& \Tor{\bu}{R\la X'_1\ra}kl
\ar@{->}[rrr]^{\Tor{\bu}{\vf\otimes_\beta\varkappa}{\ov\vf}l}
&&& \Tor{\bu}{S\la
Y'_1\ra}ll
\\
& \Tor{\bu}{A\la V\ra\la X'_1\ra}kl
\ar@{->}[u]^{\cong}
\ar@{->}[rrr]
\ar@{->}[d]^{\cong}
&&& \Tor{\bu}{B\la W\ra\la Y'_1\ra}ll
\ar@{->}[u]^{\Tor{\bu}{\breve\sigma\otimes_B B\la Y'_1\ra}ll}_{\cong}
\ar@{->}[d]_{\cong}
\\
& \Tor{\bu}{k\la V\ra}kl
\ar@{->}[rrr]_{\Tor{\bu}{\varPhi\otimes_\beta\ov\vf}{\ov\vf}l}
&&& \Tor{\bu}{l\la W\ra}ll
\\
\Tor{\bu}{k[X]}kl
\ar@{->}[rrrrr]^{\Tor{\bu}{\phi\otimes_\beta\ov\vf}{\ov\vf}l}
\ar@{->}[ur]^(0.6){\Tor{\bu}{\kappa\otimes_A k}kl}_{\cong}
&&&&& \Tor{\bu}{l[Y]}ll
\ar@{->}[ul]_(0.6){\Tor{\bu}{\lambda\otimes_B l}ll}^{\cong}
}
\]
where the non-horizontal maps are bijective by \ref{DG Tor}.
All paths commute, except possibly those around two trapezoid
with horizontal bases---the one on the floor and the larger of the
pair at the ceiling.  Composing either path from $\Tor{\bu}{A[X]\la
X'_1\ra}kl$ to $\Tor{\bu}{B\la W\ra\la Y'_1\ra}ll$ with the isomorphism
$\Tor{\bu}{\breve\sigma\otimes_BB\la Y'_1\ra}ll$ we get the same map, so
the upper trapezoid commutes.  Using this, we see that the isomorphism
$\Tor{\bu}{A[X]\otimes_A{\epsilon}}kl$, composed with either path from
$\Tor{\bu}{k[X]}kl$ to $\Tor{\bu}{l\la W\ra}ll$, yields the same map,
so the lower trapezoid commutes as well.

We inflate this trapezoid to a diagram of homomorphisms of graded algebras
\[
\xymatrixrowsep{2.0pc}
\xymatrixcolsep{1.0pc}
\xymatrix{
&& \Tor{\bu}{k\la V\ra}kl
\ar@{->}[rrr]^{\Tor{\bu}{\varPhi\otimes_\beta\ov\vf}{\ov\vf}l}
\ar@{->}[d]^{\cong}
&&&
\Tor{\bu}{l\la W\ra}ll
\ar@{->}[d]_{\cong}
\\
&& \hark l{l\la V\ra}
\ar@{->}[rrr]^{\hark l{\varPhi\otimes_\beta\ov\vf}}
&&& \hark
l{l\la W\ra}
\\
&& \hark l{l[X]}
\ar@{->}[rrr]_{\hark l{\phi\otimes_\beta\ov\vf}}
\ar@{->}[u]_{\hark
l{\kappa\otimes_A l}}
&&& \hark l{l[Y]}
\ar@{->}[u]^{\hark l{\lambda\otimes_B l}}
\\
\Tor{\bu}{k[X]}kl
\ar@{->}[rrrrrrr]^{\Tor{\bu}{\phi\otimes_\beta\ov\vf}{\ov\vf}l}
\ar@/^1pc/[uuurr]^-{\Tor{\bu}{\kappa\otimes_A k}kl}_-{\cong}
\ar@{->}[urr]_{\cong}
&&&&&&& \Tor{\bu}{l[Y]}ll
\ar@/_1pc/[uuull]_-{\Tor{\bu}{\lambda\otimes_B l}ll}^-{\cong}
\ar@{->}[ull]^{\cong}
}
\]
 From \ref{bar} we know that the maps pointing inward from the
corners are bijective, and that the upper rectangle, both triangles,
and the inner trapezoid commute.  We conclude that the lower rectangle
commutes and its vertical arrows are bijective.

Referring to \ref{bar} again, we note that all maps in the lower rectangle
are induced by morphisms of DG \G-algebras, and that $\bark l{l\la V\ra}$
and $\bark l{l\la W\ra}$ are admissible.  As $\bark l{\kappa\otimes_A l}$
and $\bark l{\lambda\otimes_B l}$ are quasiisomorphisms, it follows that
the DG \G-algebras $\bark l{l[X]}$ and $\bark l{l[Y]}$ are admissible
and the maps in the lower rectangle are isomorphisms of \G-algebras.
To finish the proof we remark, with a final reference to \ref{bar},
that the upper rectangle is formed by homomorphisms of \G-algebras.
 \end{proof}

\begin{chunk}
\label{spectral}
Let $C$ be a connected DG algebra over $l$.  The differential $\ov\dd$ of
the bar construction $\bark lC$ has the form $\dd'+\dd''$, where
\begin{gather*}
\dd'([c_1|c_2|\cdots|c_p])= \sum_{j=1}^{p-1}(-1)^{|c_1|+\cdots+|c_{j}|+j}
[c_1|\cdots|c_{j}c_{j+1}|\cdots|c_p]\\
\dd''([c_1|c_2|\cdots|c_p])= \sum_{j=1}^p(-1)^{|c_1|+\cdots+|c_{j-1}|+j}
[c_1|\cdots|\dd(c_j)|\cdots|c_p]
\end{gather*}
Thus, the $l$-span $\FF^q(C)$ of the elements $[c_1|\cdots|c_p]$ of
degree at most $(p+q)$ for $p=0,1,2\dots$, is a subcomplex of $\bark
lC$.  The page ${}^0\mbox{E}$ of the spectral sequence of the
filtration $\{\FF^q(C)\}$ is a complex of graded vector spaces with
${}^0d$ induced by $\dd'$.  It can also be obtained by tensoring with
$l$ over $C\nat$ the complex of graded $C\nat$-modules
\begin{equation}
\label{complex}
\begin{gathered}
\cdots \lra C\nat\otimes_l\Shift^p(C_{\ges1}^\natural{}^{\otimes p})
\xra{\delta_{p}}
C\nat\otimes_l\Shift^{p-1}(C_{\ges1}^\natural{}^{\otimes(p-1)})
\lra \cdots\\
\delta(\susp^p(c_0\otimes c_1\otimes\cdots\otimes c_p))=
  \sum_{j=0}^{p-1}(-1)^{j}
\susp^{p-1}(c_0\otimes\cdots\otimes c_{j}c_{j+1}\otimes\cdots\otimes c_p)
\end{gathered}
\end{equation}
where for a graded vector space $M$ we let $\Shift^p M$ denote the graded
space with $(\Shift^p M)_n=M_{n-p}$ for all $n$, and $\susp^p\col
M\to\Shift^p M$ be the degree $p$ bijection defined by the maps
$\id_{M_n}$.  The complex \eqref{complex} is the standard resolution
of $l$ by free graded $C\nat$-modules, so the spectral sequence of the
filtration $\{\FF^q(C)\}$ has
\begin{equation}
\label{sequence} \EH 1pq=\Tor p{C\nat}ll_q \implies \hark l{C}
\end{equation}
In particular, the following equalities hold:
\[
\EH 1pq =\begin{cases}
0
&\text{for }p\le0\text{ and all }q\text{ except for }(p,q)=(0,0)\,;\\
\Ind q{C}
&\text{for }p=1\text{ and all }q\,.
\end{cases}
\]
The differentials of the spectral sequence act according to the pattern
\[
\ED rpq\col\EH rpq\lra\EH r{p+r-1}{q-r}\qquad\text{for each}\qquad r\ge0
\]
so for every $q\ge1$ at the edge $p=1$ the spectral sequence defines
$l$-linear maps
\[
\hark l{C}{}_{n}\tra \EH {\infty}1{n-1}\rat \cdots \rat \EH 11{n-1}
= \Ind{n-1}{C}\,.
\]
where the kernel of the first one is the image of $\HH n{\FF^2(C)}\to
\hark l{C}{}_{n}$.  Shuffle products and divided powers in $\bark kC$
are homogeneous with respect both to degree and to weight,
cf.\ \ref{bar}, so the subspace ${\bark lC}{}^{(2)}$ of \ref{Gamma
algebra} is contained in $\FF^2(C)$.  Thus, if $\bark lC$ is
admissible, then for each $n\ge1$ the maps above define a composition
\begin{equation}
\label{analysis}
\nu^C_n\col \Gin{n}{\hark l{C}}\tra
\EH {\infty}1{n-1}\rat \cdots \rat \EH 11{n-1} = \Ind{n-1}{C}
\end{equation}
of $l$-linear homomorphisms.  Formula \eqref{naturality} yields
inclusions $\bark l{\gamma}(\FF^q(C))\subseteq\FF^q(C')$ for all $q$
and every morphism $\gamma\col C\to C'$ of connected DG algebras over
$l$. It follows that the spectral sequence \eqref{sequence} above is
natural with respect to such morphisms, and hence so are its edge
homomorphism $\nu^C_n$.
 \end{chunk}

\begin{proof}[Proof of Theorem {\em\ref{indecomposables}}]
For each $n\ge2$ we form a commutative diagram
\begin{equation*}
\xymatrixrowsep{2pc}
\xymatrixcolsep{1.3pc}
\xymatrix{
& \pie{n}{R}\otimes_k l
\ar@{=}[d]
\ar@{->}[rrrr]^-{\pie{n}{\vf}}
&&&& \pie{n}{S}
\ar@{=}[d]
\\
lX'_{n}
\ar@{->}[r]^-{\cong}
& \Gin n{\rTor{\bu}Rkl}
\ar@{->}[rrrr]^{\Gin
n{\rTor{\bu}{\vf}{\ov\vf}l}}
\ar@{->}[d]_{\cong}
&&&& \Gin n{\rTor{\bu}Sll}
\ar@{->}[d]^-{\cong}
& lY'_{n}
\ar@{->}[l]_-{\cong}
\\
& \Gin n{\hark l{l[X]}}
\ar@{->}[rrrr]^-{\Gin n{\hark l{\phi\otimes_\beta l}}}
&&&&
\Gin n{\hark l{l[Y]}}
\\
& \Ind{n-1}{l[X]})
\ar@{->}[rrrr]^-{\Ind{n-1}{\phi\otimes_\beta l}}
\ar@{<-}[u]^{\nu^{l[X]}_n}
&&&& \Ind{n-1}{l[Y]}
\ar@{<-}[u]_{\nu^{l[Y]}_n}
}
\end{equation*}
of $l$-vector spaces as follows:  The top rectangle comes from \ref{rtor}.
The middle part is obtained by stacking the commutative diagrams of
Lemmas \ref{reduced}, \ref{gammas}, and \ref{admissibles}, then taking
\G-indecomposables, as in \ref{Gamma algebra}.  The bottom rectangle
reflects the naturality of the edge homomorphisms $\nu_n$ defined in
\eqref{analysis}.

To finish the proof we show that its vertical maps are bijective.
It suffices to do this for $\nu^{l[X]}_n$.  By the isomorphisms above
and \ref{comparison}, its source and target have the same rank, so it
is enough to prove that it is surjective.  To this end we analyze the
spectral sequence \eqref{sequence}.  A well known computation, cf.\
e.g.\ \cite[(7.2.9)]{res}, gives
\[
\EH1pq=\Tor p{l[X]\nat}ll_q \cong l\la X''\ra_{p,q}
\quad\text{where}\quad
\card{X''_{m,n}}=
\begin{cases}
\card{X_n}
&\text{if } m=1\,;\\
0
& \text{otherwise}\,.
\end{cases}
\]
 From \ref{comparison} we know that $\card(X_n)=\card(X'_{n+1})$ for
$n\geq 1$, so the graded vector space associated with the bigraded
space ${}^1\mbox{E}$ is isomorphic to $l\la X'_{\ges 2}\ra$.  By Lemmas
\ref{reduced}, \ref{gammas}, and \ref{admissibles} the latter space is
isomorphic to ${\hark l{l[X]}}$.  This is the abutment of the spectral
sequence \eqref{sequence}, so it is isomorphic to the graded vector
space associated with the bigraded space ${}^\infty \mbox{E}$.  Putting
these remarks together, for each $n$ we get
\[
\sum_{p+q=n}\rank_l \EH 1pq=
\rank_l(l\la X'_{\ges 2}\ra_n)= \rank_l\HH n{\bark l{l[X]}}=
\sum_{p+q=n}\rank_l \EH {\infty}pq\,.
\]
They imply that the spectral sequence \eqref{sequence} stops on the
page ${}^1\mbox{E}$, so in the decomposition \eqref{analysis} the
injections $\EH {r+1}1{n-1}\rat\EH {r}1{n-1}$ are bijective for all
$n\ge1$ and $1\le r\le\infty$.  As a consequence, the map
$\nu_{n}^{l[X]}$ is surjective, as desired.
 \end{proof}

\section{Almost small local homomorphisms}
\label{Almost small local homomorphisms}

We introduce a class of maps of major importance for this paper.

\begin{definition}
A local homomorphism $\vf\col(R,\fm,k)\to(S,\fn,l)$ is said to be
{\em almost small\/} if the kernel of the homomorphism
$\Tor{\bu}{\vf}{\ov\vf}l\col \Tor{\bu}Rkl\to\Tor{\bu}Sll$ of graded
algebras is generated by elements of degree $1$.
 \end{definition}

The name reflects the relation of the new concept to that of {\em
small homomorphism\/}, defined in \cite{sma} by the condition that
the map $\Tor{\bu}{\vf}{\ov\vf}l$ is injective.

\begin{example}
\label{regularity}
It is proved in \cite[(4.1)]{sma} that for every ring $(R,\fm,k)$, and for
each ideal $\fa$ contained in $\fm^s$ for some sufficiently large $s$,
the canonical epimorphism $R\to R/\fa$ is small.  An effective bound on
$s$ has been found recently by Liana \c Sega.

Namely, let $G$ be the symmetric algebra of the $k$-vector space
$\fm/\fm^2$, let $\gr_\fm (R)$ be the associated graded ring of $R$,
and extend the identity map of $\fm/\fm^2$ to a homomorphism of
graded $k$-algebras $G\to\gr_\fm (R)$.  Let $\polreg R$ denote the
Castelnuovo-Mumford regularity of the graded $G$-module $\gr_\fm(R)$,
that is
\[
\polreg R=\sup_{i\in\BN}\big\{j\in\BZ\,\big|\,
\Tor iG{\gr_\fm(R)}k_{i+j}\ne0\big\}\,.
\]
By \cite[(6.2)]{Se} the epimorphism $R\to R/\fm^s$ is Golod for all
$s\geq 2+\polreg R$.  Golod homomorphisms are small by \cite[(3.5)]{sma},
so the factorization $R\to R/\fa\to R/\fm^s$ and functoriality imply
that $R\to R/\fa$ is small for every ideal $\fa$ contained in $\fm^s$.
 \end{example}

By \cite[(3.1)]{sma}, $\vf$ is small if and only if $\pies{\vf}$ is
injective.  We characterize almost smallness in similar terms,
and by means of reduced Tor-algebras, cf.\ \ref{rtor}.

\begin{proposition}
\label{small rtor}
Let $\vf\col (R,\fm,k)\to (S,\fn,l)$ be a local homomorphism.

The following conditions are equivalent.
\begin{enumerate}[\rm\quad(i)]
\item
$\vf$ is almost small.
\item
$\pie{\ges2}\vf$ is injective.
\item
$\rTor{\bu}{\vf}{\ov\vf}l$ is injective.
\end{enumerate}
\end{proposition}

\begin{proof}
Using the fact that $\Tor{\bu}{\vf}{\ov\vf}l$ is a homomorphism of {\em
Hopf $\varGamma$-algebras}, it is proved in \cite[(1.3)]{sma} that
there exists a subset $G\subset\Tor{\bu}Rkl$ such that
$\Ker(\Tor{\bu}{\vf}{\ov\vf}l)=(l\la G\ra)_{\ges1}\Tor{\bu}Rkl$ and the
following hold:
\begin{enumerate}[\rm\quad(1)]
\item
The image of $G$ in $\pies{R}\otimes_kl$ is a basis of $\Ker(\pies{\vf})$.
\item
The graded $l\la G\ra$-module $\Tor{\bu}Rkl$ is free.
\end{enumerate}
Clearly, (i) $\iff$ (ii) follows from (1).  The freeness of
$\Tor{\bu}Rkl$ over $l\la\Tor{1}Rkl\ra$, that of $\Tor{\bu}Sll$ over
$l\la\Tor{1}Sll\ra$, and (2) yield (i) $\iff$ (iii).
 \end{proof}

Vanishing of $\pie{\ges 2}R$ characterizes regularity, cf.\ Example
\ref{acyclic closure}, so we get

\begin{corollary}
If the ring $R$ is regular, then $\vf$ is almost small.  Conversely, if
the canonical surjection $\epsilon\col R\to k$ is almost small, then
$R$ is regular.
 \qed
  \end{corollary}

Using the functoriality of $\pies{\ }$, we see that the proposition
also implies

\begin{corollary}
\label{small composition}
Let $\psi\col Q\to R$ and $\vf\col R\to S$ be local homomorphisms.
\begin{enumerate}[\rm(a)]
\item
If $\psi$ and $\vf$ are almost small, then $\vf\circ\psi$ is almost small.
\item
If $\vf\circ\psi$ is almost small, then $\psi$ is almost small.
\item
If $\vf\circ\psi$ is almost small and $\pie{\ges2}\psi$ is bijective,
then $\vf$ is almost small. \qed
 \end{enumerate}
  \end{corollary}

As a further corollary, we get another example of almost small
homomorphisms.

\begin{example}
If $\vf$ is flat and $\chr(l)=2$, then $\vf$ is almost small by
Andr\'e \cite{two}.
 \end{example}

Here is what is known for flat homomorphisms in general.

\begin{remark}
\label{six}
If $\vf\col(R,\fm,k)\to(S,\fn,l)$ is a flat local homomorphism, then by
\cite[(1.1)]{cid} for every $i\ge1$ there exists an exact sequence of
$l$-vector spaces
\begin{gather*}
\xymatrixrowsep{.7pc}
\xymatrixcolsep{3pc}
\xymatrix{
0
\ar@{->}[r]
& \pie{2i}R\otimes_kl
\ar@{->}[r]^-{\pie{2i}\vf}
& \pie{2i}S
\ar@{->}[r]
& \pie{2i}{S/\fm S}
\ar@{->}[r]^-{\eth_{2i}}
&{\ }
\\
\ar@{->}[r]
& \pie{2i-1}R\otimes_kl
\ar@{->}[r]^-{\pie{2i-1}\vf}
& \pie{2i-1}S
\ar@{->}[r]
& \pie{2i-1}{S/\fm S}
\ar@{->}[r]
& 0
}
\end{gather*}
Andr\'e \cite{add} proved $\sum_{i=1}^\infty\rank_l(\eth_{2i})\le
\edim(S/\fm S)-\depth(S/\fm S)$ and conjectured that $\pie{2i-1}\vf$
is injective for all $i\ge2$.  In view of Proposition \ref{small rtor},
the conjecture can be restated to say that every flat homomorphism is
almost small.
 \end{remark}

The next proposition fails for small homomorphisms, and presents one of
the main technical reasons for working with almost small homomorphisms.

\begin{proposition}
\label{factorization}
Let $\vf\col(R,\fm,k)\to(S,\fn,l)$ be a local homomorphism.

The following maps are almost small simultaneously: $\vf$,
$\grave\vf\col R\to\wh S$, $\wh\vf\col\wh R\to\wh S$, and $\rho\col
R'\to\wh S$, where $R\xra{\dot\vf} R'\xra{\rho}\wh S$ is a regular
factorization  of $\vf$.
 \end{proposition}

\begin{proof}
The map $\grave\vf\col R\to\wh S$ is the composition of $\vf$ with
the completion map $S\to \wh S$, and also the composition of the
completion map $R\to \wh R$ with $\wh\vf$.  As $\pies{\ }$ applied
to either completion map yields an isomorphism, Corollary \ref{small
composition} shows that $\vf$, $\wh\vf$, and $\grave\vf$ are almost small
simultaneously.  Finally, $\pie{\ges 2}{R'/\fm R'}=0$ because $R'/\fm R'$
is regular, cf.\ Example \ref{acyclic closure}, so $\pie{\ges2}{\dot\vf}$
is bijective by the exact sequence of Remark \ref{six}.  Thus, $\wh\vf$
and $\grave\vf$ are almost small simultaneously by Corollary \ref{small
composition}.3.
 \end{proof}

As an application, we show how to obtain almost small homomorphisms by
factoring complete intersection homomorphisms.

\begin{corollary}
\label{small examples}
Let $\vf\col R\to(S,\fn,l)$ be a local homomorphism. If there
exists a local homomorphism $\xi\col S\to (S',\fn',l')$ such that
$\xi\circ\vf$ is c.i.\ at $\fn'$, then $\vf$ is almost small. In
particular, if $\vf$ is c.i.\ at $\fn$, then it is almost small.
  \end{corollary}

\begin{proof}
Let $R\to R'\xra{\vf'}\wh S'$ be a reduced Cohen factorization
of the composition $\wh\xi\circ\grave\vf\col R\to\wh S'$.
By hypothesis, $\Ker(\vf')$ is generated by an $R'$-regular sequence,
so \cite[(3.4.1)]{GL} shows that $\pie{n}{\vf'}$ is injective for $n=2$
and bijective for $n\ge3$. By Proposition \ref{small rtor}, the map $\vf'$
is almost small, which implies, by Proposition \ref{factorization}, that
$\xi\circ\vf$ is almost small as well.  It remains to invoke Corollary
\ref{small composition}.b.
 \end{proof}

Extending Example \ref{regularity}, we provide a numerical test for
almost smallness.

\begin{proposition}
Let $\vf\col(R,\fm,k)\to(S,\fn,l)$ be a local homomorphism.  If
\[
\length_S(S/\fn^s)=\sum_{i=0}^{s-1}\binom{e+i-1}{i}\length_R(R/\fm^{s-i})
\]
for $e=\edim(S/\fm S)$ and $s=2+\polreg R$, then $\vf$ is almost small.
  \end{proposition}

\begin{proof}
First we note that for every integer $r$ there is an inequality
\[
\length_S(S/\fn^r)\le \sum_{i=0}^{r-1}\binom{e+i-1}{i}\length_R(R/\fm^{r-i})\,.
\]
Indeed, if $R\to (R',\fm',l)\xra{\vf'}\wh S$ is a reduced Cohen
factorization, then the right hand side of the formula above is equal
to $\length_{R'}(R'/{\fm'}^r)$.  The surjective homomorphism $\vf'_r\col
R'/{\fm'}^r\to S/\fn^r$ induced by $R'\to\wh S$ yields the inequality
above, and $\vf'_r$ is bijective if and only if equality holds, that is,
if and only if $\Ker(\vf')\subseteq\fm'{}^r$.

By the preceding argument, our hypothesis implies $\Ker(\vf')\subseteq
\fm'{}^s$ with $s=2+\polreg R$. On the other hand, the associated
graded rings of $R'$ and $R$ are linked by an isomorphism of graded
$l$-algebras
\[
\gr_{\fm'}(R')\cong l\otimes_k \gr_{\fm}(R)[x_1,\dots,x_e]
\]
where $x_1,\dots,x_e$ are indeterminates.  It follows that $\polreg(R')
=\polreg(R)$.  Example \ref{regularity} now shows that the homomorphism
$\vf'$ is small.  In view of Proposition \ref{factorization}, it follows
that the homomorphism $\vf$ is almost small.
 \end{proof}

The next result is a structure theorem for morphisms of minimal models
over certain almost small homomorphisms.  A key ingredient of the proof
is the general result on the map $\pies{\vf}$ established in Theorem
\ref{indecomposables}.

\begin{theorem}
\label{egg}
Let $\rho\col(A,\fp,k)\to R$ and $\vf\col(R,\fm,k)\to S$ be surjective
homomorphisms of local rings, with $A$ regular and $\Ker(\rho)\subseteq
\fp^2$.  If $\vf$ is almost small, then there exists a commutative
diagram of morphisms of DG algebras
\[
\xymatrixrowsep{2pc}
\xymatrixcolsep{2.5pc}
\xymatrix{
\ A_{\mathstrut}\
\ar@{^{(}->}[r]
\ar@/^2pc/[rr]^{\beta}
\ar@{^{(}->}[d]
\ar@/_2pc/[dd]_{\rho}
&\,A[U]_{\mathstrut}\,
\ar@{->>}[r]^{\wt\beta}_\simeq
\ar@{^{(}->}[d]
& B_{\mathstrut}
\ar@/^2pc/[dd]^{\sigma}
\ar@{_{(}->}[d]
\\
\,A[X]\,
\ar@{^{(}->}[r]
\ar@{->>}[d]_{\wt\rho}^{\simeq}
& A[X,U,T]
\ar@{->>}[r]^{\chi}_\simeq
\ar@{->>}[d]_{\upsilon}^\simeq
& B[X,T]
\ar@{->>}[d]^{\wt\sigma}_\simeq
\\
\,R\,
\ar@{^{(}->}[r]
\ar@/_2pc/[rr]_{\vf}
& R[U,T]
\ar@{->>}[r]^{\wt\varphi}_\simeq
&
S
}
\]
where $(B,\fq,k)$ is a regular local ring, $\beta$ and $\sigma$ are
surjective homomorphisms, $\Ker(\sigma)$ is contained in $\fq^2$,
$U=U_1$, the external rows and columns are minimal models, $\chi$ and
$\upsilon$ are surjective quasiisomorphisms.
 \end{theorem}

\begin{proof}
Choose a subset $\bsa\subset A$ mapping to a basis of
$(\Ker(\vf)+\fm^2)/\fm^2$.

As $\bsa$ is part of a regular system of parameters for $A$, the local
ring $(B,\fq,k)=(A/(\bsa),\fp/(\bsa),k)$ is regular.  Since $\bsa$
is contained in the kernel of $\vf\circ\rho$, this map factors as a
composition of surjective homomorphisms $\beta\col A\to B$ and $\sigma\col
B\to S$.  The choice of $\bsa$ ensures that $\Ker(\sigma)$ is contained
in $\fq^2$.

Using \ref{models}, we form a commutative diagram of morphisms of DG
algebras
\[
\xymatrixrowsep{2pc}
\xymatrixcolsep{2pc}
\xymatrix{
A_{\mathstrut}
\ar@{->}[r]^{\beta}
\ar@/_2pc/[dd]_{\rho}
\ar@{^{(}->}[d]
&
B_{\mathstrut}
\ar@/^2pc/[dd]^{\sigma}
\ar@{_{(}->}[d]
\\
A[X]
\ar@{->}[r]^{\phi}
\ar@{->>}[d]_{\wt\rho}^{\simeq}
& B[Y]
\ar@{->>}[d]^{\wt\sigma}_{\simeq}
\\
R
\ar@{->}[r]^{\vf}
& S
}
\]
It induces morphisms of DG algebras
\begin{gather*}
\phi'\col B[X]=A[X]\otimes_AB\xra{\phi\otimes_AB}B[Y]\otimes_AB=B[Y]\,;\\
\ov\phi\col k[X]=A[X]\otimes_Ak\xra{\phi\otimes_Ak}B[Y]\otimes_Ak=k[Y]\,.
\end{gather*}
As $\vf$ is almost small, Theorem \ref{indecomposables} shows that the
$k$-linear map
\[
\ind{\ov\phi}\col\ind{k[X]}\lra\ind{k[Y]}
\]
is injective.  Thus, the set $\ind{\ov\phi}(X)$ is linearly independent
in $\ind{k[Y]}$.

Choose a subset $T$ in $B[Y]$ whose image in $\ind{k[Y]}$ extends
$\ind{\ov\phi}(X)$ to a basis.  It follows that $T$ is a set of free
variables over $k[X]$, and so the map $\ov\phi$ is injective.  This map
is equal to $\phi'\otimes_Bk$, and $\phi'$ is a map of graded free
$B$-modules, so we conclude by Nakayama's Lemma that $\phi'$ is injective
and $T$ is a set of free variables generating $B[Y]$ over $\phi'(B[X])$.
Changing variables in $B[Y]$, we replace $B[Y]$ by $B[X,T]$ and $\phi'$
by the canonical inclusion $B[X]\hra B[X,T]$.

Let $A[U]$ be the Koszul complex with $U=U_1$ and $\dd(U)=\bsa$, and
set $A[X,U]=A[X]\otimes_AA[U]$.  These algebras appear in a commutative
diagram of DG algebras
\[
\xymatrixrowsep{2pc}
\xymatrixcolsep{2.5pc}
\xymatrix{
A_{\mathstrut}
\ar@{^{(}->}[r]
\ar@{^{(}->}[d]
& A[U]_{\mathstrut}
\ar@{^{(}->}[d]
\ar@{->>}[r]^{\wt\beta}_{\simeq}
&
B_{\mathstrut}
\ar@{_{(}->}[d]
\\
\,A[X]\,
\ar@{^{(}->}[r]
\ar@{=}[d]
& A[X,U]_{\mathstrut}
\ar@{^{(}->}[d]
\ar@{->>}[r]^{\chi^0}_{\simeq}
& B[X]_{\mathstrut}
\ar@{_{(}->}[d]
\\
\,A[X]\,
\ar@{^{(}->}[r]
& A[X,U,T]
\ar@{->>}[r]^{\chi}_{\simeq}
& B[X,T]
}
\]
where $\wt\beta$ is the canonical augmentation and $\chi^0=
A[X]\otimes_A\wt\beta$.  Since $\bsa$ is an $A$-regular sequence,
$\wt\beta$ is a quasiisomorphism; by \ref{quisms}, $\chi^0$ is a
quasiisomorphism as well.

The map $\chi$ is built inductively, starting with $\chi^0$.  Using the
inclusions $\dd(T_1) \subseteq\Ker(\sigma)\subseteq\fq^2$, pick for
each $t\in T_1$ an element $p_t\in\fp^2$ with $\beta(p_t)=\dd(t)$. Let
$A[T_1]$ be the Koszul complex with $\dd(t)=p_t$ and $\chi^1$ the
morphism of DG algebras
\[
A[X,U,T_1]=A[X,U]\otimes_AA[T_1] \xra{\chi^0\otimes_AA[T_1]}
B[X]\otimes_AA[T_1]=B[X,T_1]\,.
\]
By \ref{quisms}, $\chi^1$ is a surjective quasiisomorphism.  Assume
next that a surjective quasiisomorphism $\chi^n\col A[X,U,T_{\les
n}]\to B[X,T_{\les n}]$ is available for some $n\ge1$.  For each
$t\in T_{n+1}$ we pick a cycle $z_t\in A[X,U,T_{\les n}]_n$ such that
$\chi^n(z_t)=\dd(t)$, then we set $A[X,U,T_{\les n+1}]=A[X,U,T_{\les
n}][T_{n+1}\var\dd(t)= z_t]$ and define $\chi^{n+1}$ to be the extension
of $\chi^n$ satisfying $\chi^{n+1}(t)=t$ for all $t\in T_{n+1}$.
It is easy to verify that this map is a surjective quasiisomorphism,
cf.\ also \cite[(7.2.10)]{res}.  Taking direct limits, we obtain the
surjective quasiisomorphism $\chi$ displayed in the diagram.

The diagram above provides the two upper squares of the diagram in
the theorem.  Its lower left square is obtained by base change along
$\wt\rho$.  For its lower right square, we factor $\wt\sigma\circ\chi$
through $\upsilon= \wt\rho\otimes_{A[X]}A[X,U,T]$ to get a surjection
$\wt\vf\col R[U,T]\to S$.

The top row and two side columns of the diagram are  minimal models
by construction.  Furthermore, $\chi$, $\upsilon$, and $\wt\vf$ are
surjective quasiisomorphisms: the first by construction, the second by
\ref{quisms}, and the third due to the commutativity of the diagram.
Since the differential of $R[U,T]$ is induced by that of $A[X,U,T]$,
to prove that the lower row is a minimal model it suffices to establish
that the differential on $A[X,U,T]$ is decomposable.

For any $y\in A[X,U,T]$ with $|y|=n+1\ge2$, write $\dd(y)$ in the form
\[
\dd(y)=\sum_{x\in X_n}c_x x+\sum_{u\in U_n}b_u u+
\sum_{t\in T_n}a_t t+w \in A[X,U,T]
\]
with $a_t, b_u, c_x\in A$ and $w\in(X,U,T)^2A[X,U,T]$.  In the
resulting equality
\[
\dd(\chi(y))=\sum_{x\in X_n}\beta(c_x)x+\sum_{t\in T_n}\beta(a_t)t
+\chi(w)\in B[X,T]
\]
we have $\chi(w)\in(X,T)^2B[X,T]$.  The differential of $B[X,T]$ is
decomposable, so for all $t\in T_n$ and $x\in X_n$ we obtain
$\beta(a_t),\beta(c_x)\in\fq$, that is, $a_t,c_x\in\fp$. Since $U=U_1$,
we have $b_u=0$ unless $n=1$.  If $n=1$, then $w=0$, so the equality
$\dd^2(y)=0$ yields
\[
\sum_{u\in U_1}b_u\dd(u)=
-\sum_{x\in X_1}c_x\dd(x)-\sum_{t\in T_1}a_t\dd(t)\,.
\]
By construction, we have $\dd(x)\in\Ker(\rho)\subseteq\fp^2$ for all
$x\in X_1$ and $\dd(t)=p_t\in\fp^2$ for all $t\in T_1$, so the last
equality yields $\sum_{u\in U_1}b_u\dd(u)\in\fp^2$.  As $\dd(U_1)=\bsa$
is part of a regular system of parameters, this implies $b_u\in \fp$
for all $u\in U_1$, so the differential of $A[X,U,T]$ is decomposable.
 \end{proof}

\section{Weak category of a local homomorphism}
\label{Weak category of a local homomorphism}

We introduce a notion motivated by F\'elix and Halperin's \cite[(4.3)]{FH}
definition of rational Lusternik-Schnirelmann category $\mbox{cat}_0(X)$
of a simply connected CW complex $X$ of finite type.  Weak category
captures a Looking Glass \cite{AH1} image of a corollary of the Mapping
Theorem:  If $\mbox{cat}_0(X)\le s$, then by \cite[(5.1)]{FH} for each
$n\ge2$ the $n$-connected cover $X_n$ of $X$ satisfies $\mbox{cat}_0(X_n)
\le s$, hence by \cite[(4.10)]{FH} the product of any $(s+1)$ cohomology
classes in $\HH{\ges1}{X_n;\BQ}$ is equal to $0$.

\begin{definition}
\label{algebra cat}
If $(B,\fq,k)$ is a local ring and $B[V]$ is a semifree extension with
decomposable differential, then we define a notion of {\em weak
category\/} by the formula
\[
\wcat(B[V])=\inf\left\{s\in\BN \left\vert\,\begin{gathered}
\text{for each $n\geq 2$ the product of any}\\
\text{$(s+1)$ elements of positive degree}\\
\text{in $\hh{B[V]/(\fq,V_{\sles n})}$ is equal to $0$}
\end{gathered}\,\right.\right\}.
\]
 \end{definition}

Finite weak category can often be detected by using a variant of
\cite[(1.2)]{hh}:

\begin{proposition}
\label{mapping theorem}
If $(B,\fq,k)$ is a local ring, $B[V]$ a DG algebra with decomposable
differential, and $B[V]\to S[W]$ a surjective morphism of DG algebras,
then
\[
\wcat(S[W])\le\sup\{s\in\BN\mid\HH{s}{B[V]/\fq B[V]}\ne0\}\,.
\]
 \end{proposition}

\begin{proof}
Set $k[V]=B[V]/\fq B[V]$ and $k[W]=S[W]/\fq S[W]$.  All DG algebras
under consideration are images of $B[V]$, so their differentials are
decomposable.  We have $\wcat(S[W])=\wcat(k[W])$ by definition.  As the
induced morphism $k[V]\to k[W]$ is surjective, we get $\wcat(k[W])\le
\sup\{s\in\BN\mid\HH{s}{k[V]}\ne0\}$ from \cite[(1.2)]{hh}.
 \end{proof}

\begin{definition}
\label{hom cat}
If $\vf\col R\to S$ is a local homomorphism and $R\to R'[U]\to \wh S$
is a minimal model of $\grave\vf$, then we define the {\em weak
category} of $\vf$ by the equality
\[
\wcat(\vf)= \wcat(R'[U])\,.
\]
Proposition \ref{independence} shows that it does not depend on the
choice of minimal model.
 \end{definition}

By \cite[\S 3]{lci}, the sequence $(\ve_n(\vf))$ is positive and grows
exponentially when $\vf$ is not c.i.\ and $\fd_RS$ is finite.  A close
reading of the proofs shows the last condition is used only to ensure
$\wcat(\vf)<\infty$, cf.\ Theorem \ref{ffd}, so at no further expense
we get

\begin{theorem}
\label{growth}
Let $\vf\col R\to (S,\fn,l)$ be a local homomorphism.

If $\wcat(\vf)$ is finite, then the following conditions are equivalent.
\begin{enumerate}[\rm\quad(i)]
\item
$\vf$ is not complete intersection at $\fn$\,.
\item
$\ve_n(\vf)>0$ for all $n\geq 2$\,.
\item
$\limsup_n\root{n}\of{\ve_n(\vf)}>1$\,.
\item
There exist a real number $c > 1$ and a sequence of integers $s_j$ with
\begin{gather*}
0<2s_j\le s_{j+1}\le(\wcat(\vf)+1)s_j
\quad\text{and}\quad \ve_{s_j}(\vf) > c^{s_j}
\quad\text{for all}\quad j \ge 1\,.
\end{gather*}
\end{enumerate}
 \end{theorem}

\begin{proof}
Proposition \ref{ci map} shows that (ii) or (iii) implies (i).

If (i) holds, then $\ve_n(\vf)\ne0$ for $n=2,3$ by Propositions
\ref{regular map} and \ref{ci map}, and $\ve_n(\vf)\ne0$ $n\ge4$ by the
{\em proof\/} of \cite[(3.4)]{lci}; thus, (i) implies (ii).

The {\em proof\/} of \cite[(3.10)]{lci} shows that (i) implies (iv).
 \end{proof}

\begin{corollary}
\label{rigidity}
The following conditions are equivalent.
\begin{enumerate}[\rm\quad(i)]
\item
$\vf$ is complete intersection at $\fn$\,.
\item
$\wcat(\vf)<\infty$ and $\ve_n(\vf)=0$ for some $n\geq 2$\,.
\item
$\wcat(\vf)<\infty$ and $\limsup_n\root{n}\of{\ve_n(\vf)}\le1$\,.
\item
$\wcat(\vf)=0$\,.
\end{enumerate}
 \end{corollary}

\begin{proof}
The theorem shows that (i) follows from either (ii) or (iii).

Let $R\to R'[U]\to \wh S$ be a minimal model of $\grave\vf\col R\to\wh S$.
Proposition \ref{ci map} shows that $\vf$ is c.i.\ at $\fn$ if and only
if $U=U_1$.  Thus, (i) implies (iv) by definition of $\wcat(\vf)$.
Conversely, if (iv) holds, then $U=U_1$ by \ref{decomposable}, hence
(i) holds.

Finally, conditions (i) and (iv) imply (ii) and (iii) by Proposition
\ref{ci map}.
 \end{proof}

Next we establish a most important property of almost small homomorphisms.

\begin{theorem}\label{finite wcat}
If a local homomorphism $\vf\col R\to S$ is almost small, then
\[
\wcat(\vf)\leq \edim S - \depth S\,.
\]
 \end{theorem}

\begin{proof}
Let $R\xra{\dot\vf} R'\xra{\vf'}\wh S$ be a Cohen presentation
of $\grave\vf$.  As $R'\xra{=}R'\xra{\vf'}\wh S$ is a Cohen
presentation of $\grave\vf'=\vf'$, we have $\wcat(\vf)=\wcat(\vf')$.
On the other hand, the map $\vf'$ is almost small by Proposition
\ref{factorization}. Furthermore, $\edim \wh S=\edim S$ and $\depth\wh
S=\depth S$.  Thus, we may assume that $R$ and $S$ are complete and that
the local homomorphism $\vf\col (R,\fm,k) \to (S,\fn,k)$ is surjective
and almost small.

Choose a regular local ring $(A,\fp,k)$ and a surjective homomorphism
$\rho\col A\to R$ with $\Ker(\rho)\subseteq\fp^2$.  From Theorem
\ref{egg} we get a minimal model $R\to R[U,T]\to S$ where $U=U_1$,
together with a minimal model $B\to B[X,T]\xra{\wt\sigma}S$ where
$(B,\fq,k)$ is a regular local ring and $\Ker(\wt\sigma_0)\subseteq
\fq^2$, linked for each $n\ge2$ by isomorphisms
\[
\frac{R[U,T]}{(\fm,U,T_{<n})R[U,T]}\cong
\frac{A[X,U,T]}{(\fp,X,U,T_{<n})A[X,U,T]}\cong
\frac{B[X,T]}{(\fq,X,T_{<n})B[X,T]}\,.
\]
We also have $\HH i{B[X,T]/\fq B[X,T]}=\Tor iBSk$ by definition, and $\Tor
iBSk=0$ for $i>\dim B - \depth S$ by the Auslander-Buchsbaum Equality.
Thus, Proposition \ref{mapping theorem} yields $\wcat(\vf)= \wcat(R[U,T])
\le\dim B-\depth S$.  It remains to note that $\dim B=\edim S$ because
$\Ker(\wt{\sigma}_0)$ is contained in $\fq^2$.
 \end{proof}

For completeness, we deduce \cite[(3.8)]{lci} from Proposition
\ref{mapping theorem}.

\begin{theorem}
\label{ffd} If $\vf\col(R,\fm,k)\to S$ is a local homomorphism, then
\[
\wcat(\vf)\leq\edim(S/\fm S)+\fd_RS\,.
\]
 \end{theorem}

\begin{proof}
There is nothing to prove unless $f=\fd_RS$ is finite.  Let $R\to
R'[U]\to\wh S$ be a minimal model of $\grave\vf$ with $\edim R'/\fm
R'=\edim S/\fm S$; call this number $e$.  By Proposition \ref{mapping
theorem}, it suffices to show $\HH i{R'[U]/\fm'R'[U]}=0$ for $i>e+f$.

Since $R'[U]$ is a flat resolution of $\wh S$ over $R$, we have $\HH
i{F}=\Tor iRkS=0$ for $i>f$ and $F=R'[U]/\fm R'[U]$.  Now note that $F$
is a complex of free modules over the regular local ring $\ov R =R'/\fm
R'$, and that $\dim\ov R=e$ by the minimality of the Cohen
factorization. For $i=0,\dots,e-1$ form exact sequences of complexes
\[
0\lra F/(x_1,\dots,x_i)\xra{x_{i+1}}F/(x_1,\dots,x_i)\lra
F/(x_1,\dots,x_{i+1})\lra0
\]
where $x_1,\dots,x_e$ is a regular system of parameters of $\ov R$.
 From their homology exact sequences, one sees by induction on $i$
that the homology of the complex $F/(x_1,\dots,x_e)=R'[U]/\fm'R'[U]$
vanishes in degrees greater than $e+f$.
 \end{proof}

\section{Andr\'e-Quillen homology of local homomorphisms}
\label{Andre-Quillen homology}

In this section we prove local versions of our results on Andr\'e-Quillen
homology.  When using the general properties of the theory we take
Andr\'e's monograph \cite{hca} as standard reference.  In addition,
we heavily draw on some results from \cite{lci}, {\em verbatim\/} or
in variants.  We recall them below.

\begin{chunk}
\label{low dimensions}
For each local homomorphism $\psi\col Q\to (R,\fm,k)$, by
\cite[(4.3)]{lci} one has
\[
\rank_k\aq nRQk=\ve_{n+1}(\psi) \text{ for }
\begin{cases}
2\le n < \infty
&\text{if}\quad\chr k=0\,;\\
2\le n\le 2p-1
&\text{if}\quad\chr k=p>0\,.
\end{cases}
\]
 \end{chunk}

\begin{chunk}
\label{flat ci}
For each local homomorphism $\psi\col Q\to (R,\fm,k)$ the following are
equivalent.
\begin{enumerate}[{\quad\rm(i)}]
\item
$\aq nRQk=0$ for all $n\gg0$, and $\fd_QR<\infty$\,.
\item
$\aq nRQk=0$ for all $n\geq 2$\,.
\item
$\aq 2RQk=0$\,.
\item
$\aq nRQk=0$ for some $n\ge2$ such that $\lfloor \frac n2\rfloor!\ne0\in
k$, and $\fd_QR<\infty$\,.
\item
$\psi$ is complete intersection at $\fm$\,.
\end{enumerate}
Indeed, (ii), (iii), and (v) are equivalent by \cite[(1.8)]{lci}.  If
$\psi$ is c.i.\ at $\fm$, then $\fd_QR$ is finite, cf.\ \cite[(3.2)]{afh},
so (ii) and (v) imply (i) and (iv). Conversely, (i) implies (v) by
\cite[(4.4)]{lci}, while (iv) implies (v) by \cite[(3.4)]{lci} via the
equalities in \ref{low dimensions}.
 \end{chunk}

\begin{chunk}
\label{regular factorization}
Let $\xi\col(R,\fm,k)\to(R^*,\fm^*,k^*)$ be a flat local homomorphism
such that $\fm^*=\fm R^*$.  If in a commutative square of local
homomorphisms
\[
\xymatrixrowsep{2pc}
\xymatrixcolsep{3pc}
\xymatrix{
&P
\ar@{->}[dr]
\\
Q
\ar@{->}[r]^-{\psi}
\ar@{->}[ur]^{\iota}
& R
\ar@{->}[r]^-{\xi}
& R^*
}
\]
the upper path is a regular factorization of the composition
$\xi\circ\psi$, then the canonical maps $\aq n{\xi}{\iota}{k^*}\col \aq
nRQ{k^*}\to\aq n{R^*}P{k^*}$ are bijective for all $n\ge2$.

Indeed, the argument for \cite[(1.7)]{lci} carries over with only
notational changes.
 \end{chunk}

We now present our main local result, describing c.i.\ homomorphisms
in terms parallel to those in \ref{flat ci}, but {\em without\/}
the hypothesis of finite flat dimension.

\begin{theorem}
\label{almost quillen}
If $\psi\col Q\to (R,\fm,k)$ is a local homomorphism, then the
following conditions are equivalent.
\begin{enumerate}[{\quad\rm(i)}]
\item
$\aq nRQk=0$ for all $n\gg 0$, and $\psi$ is almost small.
\item
$\aq nRQk=0$ for all $n\gg 0$ and $\wcat(\psi)<\infty$.
\item
$\aq nRQk=0$ for some $n\ge2$ with $\lfloor \frac n2\rfloor!\ne0\in k$,
and $\psi$ is almost small\,.
\item
$\aq nRQk=0$ for some $n\ge2$ with $\lfloor \frac n2\rfloor!\ne0\in k$,
and $\wcat(\psi)<\infty$\,.
\item
$\psi$ is complete intersection at $\fm$\,.
\end{enumerate}
 \end{theorem}

\begin{proof}
If $\psi$ is c.i.\ at $\fm$, then $\aq nRQk=0$ for all $n\ge2$ by
\ref{flat ci}, and $\psi$ is almost small by Corollary \ref{small
examples}, so (v) implies (i) and (iii).  If $\psi$ is almost small, then
it has finite weak category by Theorem \ref{finite wcat}, so (i) implies
(ii) and (iii) implies (iv).  If (ii) holds, then the {\em proof\/} of
\cite[(4.4)]{lci} shows that $\limsup_n\root n\of{\ve_n(\psi)}\le1$, so
(v) holds by Corollary \ref{rigidity}.  In view of \ref{low dimensions},
the same corollary shows that (iv) implies (v).
 \end{proof}

To continue, we recall some general computations of Andr\'e-Quillen
homology.

\begin{chunk}
Let $\vf\col R\to S$ be a homomorphism of commutative rings, and $N$
an $S$-module.

Andr\'e \cite[(16.1)]{simp} constructs a universal coefficients spectral
sequence
\[
{}^2\operatorname{E}_{p,q}= \Tor pS{\aq qSRS}N\implies\aq {p+q}SRN\,.
\]
Thus, if for some $m\in\BZ$ the $S$-module $\aq nSRS$ is flat for
all $n\le m$, then
\begin{equation}
\label{universal}
\aq nSRN\cong\aq nSRS\otimes_SN \qquad\text{for all}\qquad n\le m+1\,.
\end{equation}

Set $\fa=\Ker(\vf)$. There are isomorphisms of $S$-modules
\begin{gather}
\label{deg1}
\aq 1SRN\cong\fa/\fa^2\otimes_RN\cong\Tor1RSN\,;\\
\label{deg2}
\aq 2SRS\cong\frac{\Tor2RSS}{\Tor1RSS\cdot\Tor1RSS}\,,
\end{gather}
where the second one is elementary, and the other two come from
\cite[(6.1), (15.8)]{hca}.
 \end{chunk}

We give ``concrete'' descriptions of split c.i.\ local homomorphisms.

\begin{chunk}
\label{factor}
Let $\vf\col(R,\fm,k)\to (S,\fn,k)$ be a local homomorphism, set
$\Ker(\vf)= \fa$, and let asterisks ${\,}^*$ denote $\fa$-adic completion.
If $\psi\col S\to R$ is a section of $\vf$, then there exist a set $\bsx$
of formal indeterminates over $S$ and a commutative diagram
\[
\xymatrixrowsep{2pc}
\xymatrixcolsep{2.5pc}
\xymatrix{
&{\ }R^*
\ar@{->}[dr]^{\varphi^*}
\\
S
\ar@{->}[r]^{\iota}
\ar@{->}[ur]^{\psi^*}
& {S[[\bsx]]}/(\bsf)
\ar@{->}[u]^{\cong}
\ar@{->}[r]^{\pi}
& S
}
\]
of homomorphisms of rings, where the isomorphism is induced by a
surjective map $\rho\col S[[\bsx]]\to R^*$ with $\Ker(\rho)\subseteq
\fn(\bsx) +(\bsx)^2$, $\bsf$ is a minimal set of generators of
$\Ker(\rho)$, $\pi$ is the surjection with $\Ker(\pi)=(\bsx)$, and $\iota$
is the natural injection.

Indeed, as $S$ is discrete in the $\fa$-adic topology, completion
yields homomorphisms of rings $\psi^*\col S\to R^*$ and $\varphi^*\col
R^*\to S$ whose composition is the identity map of $S$.  Pick a
set $\bsa=\{a_1,\dots,a_e\}\subset R^*$ minimally generating
$\Ker\varphi^*$, let $\bsx=\{x_1,\dots,x_e\}$ be a set of formal
indeterminates, and let $\rho\col S[[\bsx]]\to R^*$ be the unique
homomorphism of $S$-algebras mapping $x_i$ to $a_i$ for each $i$.
It is necessarily surjective, and in view of the choice of $\bsa$ we
have $\Ker(\rho)\subseteq(\fn+(\bsx))^2\cap(\bsx)=\fn(\bsx)+(\bsx)^2$.
 \end{chunk}

\begin{proposition}
\label{split ci}
In the notation of {\rm\ref{factor}} the following hold.
\begin{enumerate}[\rm(a)]
\item
$\bsf$ is a regular sequence if and only if $\psi$ is c.i.\ at $\fm$.
\item
$\bsf$ is a regular sequence contained in $(\bsx)^2$ if and only if
the $S$-modules $\aq 1SRS$ and $\aq 2SRS$ are free, and $\aq 3SRS=0$.
\item
$\bsf=\varnothing$ if and only if the $S$-module $\aq 1SRS$ is free,
$\aq 2SRS=0$, and $\aq 3SRS=0$.
 \end{enumerate}
  \end{proposition}

\begin{proof}
Set $P=S[[\bsx]]$.  The ring $R^*$ is local with maximal ideal $\fm R^*$
and residue field $R^*/\fm R^*\cong k$, and $R^*$ is flat over $R$.  Thus,
\cite[(4.54)]{hca} yields $\aq nSRS\cong\aq nS{R^*}S$ for all $n\in\BZ$,
so for the rest of the proof we may assume $R=R^*$.

As $\aq 2RSk\cong \aq 2{R}Pk$ by \ref{regular factorization}, (a)
follows from \ref{flat ci}.

Since $\bsx$ is $P$-regular, $\aq nSP-=0$ for $n\ge2$, so the
Jacobi-Zariski exact sequence \cite[(5.1)]{hca} generated by $P\to R\to S$
yields isomorphisms of functors
\begin{equation}
\label{degree shift} \aq nSR-\cong\aq {n-1}RP- \qquad\text{ for }\qquad n\ge3
\end{equation}
and, in view of the isomorphism \eqref{deg1}, also an exact
sequence
\begin{equation}
\label{four terms}
0\lra\aq2SRS\lra({(\bsf)}/{(\bsf)^2})\otimes_{R}S\xra{\delta}
{(\bsx)}/{(\bsx)^2}\lra\aq1SRS\lra0
\end{equation}
of $S$-modules.  Note that the $S$-module $(\bsx)/(\bsx)^2$ is free,
and
\begin{equation}
\label{inclusion}
\im(\delta)=((\bsf)+(\bsx)^2)/(\bsx)^2\subseteq
(\fn(\bsx)+(\bsx)^2)/(\bsx)^2=\fn\big((\bsx)/(\bsx)^2\big)\,.
\end{equation}

If $\bsf$ is a regular sequence, then $\aq 2RPS=0$ and $(\bsf)/(\bsf)^2$
is free over $R$, so the ``only if'' parts of (b) and (c) are now clear.
To obtain converses, we assume that $\aq nSRS$ is free for $n\le2$,
and $\aq 3SRS=0$.  From \eqref{degree shift} and \eqref{universal} we get
\[
\aq 2RPk\cong\aq 3SRk\cong\aq 3SRS\otimes_Sk=0\,,
\]
from where we conclude that $\bsf$ is a regular sequence, cf.\ \ref{flat
ci}.  Due to \eqref{four terms} and \eqref{inclusion}, the freeness
of $\aq 1SRS$ implies $\delta=0$, that is, $\bsf \subseteq(\bsx)^2$
and $\aq2SRS\cong(\bsf)/(\bsf)^2$.  It is now clear that the ``if''
parts of (b) and (c) hold as well.
 \end{proof}

The next theorem is a local version of Theorem II.  The proof draws on
results obtained above and on earlier results from \cite{hh}.  One of
them is for a homomorphism $\vf\col (R,\fm,k) \to (S,\fn,k)$ of local
rings that is {\em large\/} in the sense of Levin \cite{Le2}, meaning
that the map $\Tor{\bu}{\vf}{\ov\vf}k\col \Tor{\bu}Rkk\to\Tor{\bu}Skk$
is surjective.

\begin{theorem}
\label{local second}
For an algebra retract $S\xra{\psi}(R,\fm,k)\xra{\vf}S$ of local rings
the following conditions are equivalent.
\begin{enumerate}[{\quad\rm(i)}]
\item
The $S$-algebra $\Tor{\bu}RSS$ is finitely generated.
\item
For every $S$-algebra $T$ there is an isomorphism of graded
$T$-algebras
\[
\Tor{\bu}RS{T}\cong
\big({\ts\bigwedge}_{S}D_1\otimes_S\sym_{S}D_2\big)\otimes_ST
\]
where $D_1$ and $D_2$ are free $S$-modules concentrated in degrees $1$
and $2$, respectively. Moreover, if $\chr(k)>0$, then $D_2=0$.
\item
The $S$-modules $\aq 1SRS$ and $\aq 2SRS$ are free, and $\aq
3SRS=0$. Moreover, if $\chr(k)>0$, then $\aq 2SRS=0$.
\item
The sequence $\bsf$ in the commutative diagram constructed in
{\em\ref{factor}} is contained in $(\bsx)^2$. Moreover, if $\chr(k)>0$,
then $\bsf=\emptyset$.
 \end{enumerate}
  \end{theorem}

\begin{proof}
Suppose that $\chr(k)>0$.  Conditions (iii) and (iv) are then equivalent
by Proposition \ref{split ci}.c.  Functoriality in the ring argument
shows that $\vf$ is large.  For every large homomorphism, conditions (i),
(ii) and (iv) are equivalent by \cite[(3.1)]{hh}.

Now consider the case when $\chr(k)=0$.  Conditions(iii) and (iv) are
equivalent by Proposition \ref{split ci}.b.  The equivalence of (i),
(ii), and (iv) is established in \cite[(4.1)]{hh} under an additional
hypothesis, $R$ has finite flat dimension over $S$, which is used only
once, to conclude that $\aq n{R^*}{S[[\bsx]]}k=0$ for all $n\gg0$ implies
that $\bsf$ is regular (cf.\ \cite[p.\ 163]{hh}).  By Theorem \ref{almost
quillen}, the same conclusion holds when the map $\rho\col S[[\bsx]]\to
R^*$ is almost small.  To see that it is, note that the composition
$S[[\bsx]]\xra{\rho} R^*\xra{\psi^*}S$ is complete intersection at $\fn$,
and apply Corollary \ref{small examples}.
 \end{proof}

The proof above raises the question whether the first three conditions
are equivalent for large homomorphisms in characteristic $0$.  Here is
what we know.

\begin{chunk}
\label{large}
Let $\vf\col R\to (S,\fn,k)$ be a large homomorphism, where $\chr(k)=0$.
If $R\to R[U]\to S$ is a minimal model of $\vf$, then
results of Quillen, \cite[(9.5), (10.1)]{Qu}, imply
\[
\sup\{n\in\BN\mid \aq nSRk\ne0\}=
\sup\{n\in\BN\mid U_n\ne\varnothing\}\,.
\]
Thus, \cite[(3.2)]{hh} shows that if $\Tor{\bu}RSS$ is finitely generated,
then $\aqdim RS$ is finite, and \cite[(3.1)]{acy} proves that in this
case $\aqdim RS$ is equal to $1$ or is even.
 \end{chunk}

\section{Andr\'e-Quillen dimension}
\label{Andre-Quillen dimension}

In this section we prove the theorems stated in the Introduction.
The local case was essentially settled in Section \ref{Andre-Quillen
homology}, but reduction to that case needs some attention, as weak
category and almost small homomorphism are intrinsically local notions.

We start by recording a slight extension of a result of Andr\'e
\cite{hca}.  For each $\fn\in\Spec S$, we let $k(\fn)$ denote the residue
field $S_\fn/\fn S_\fn$.

\begin{chunk}
\label{localization}
If $Q\to R\to S$ are homomorphisms of noetherian rings and $j$ is
a non-negative integer, then the following conditions are equivalent.
\begin{enumerate}[{\quad\rm(i)}]
\item
$\aq nRQ-=0$ on the category of $S$-modules for all $n\ge j$.
\item
$\aq n {R_{\fm}}{Q_{\fn\cap Q}}{k(\fm)}=0$ for all $n\ge j$, all
$\fn\in\Spec S$, and $\fm=\fn\cap R$.
 \end{enumerate}

Indeed, \cite[(4.58), (5.27)]{hca} produce for each $n$ canonical
isomorphisms
\[
\aq nRQ{k(\fn)}\cong
\aq n{R_{\fm}}Q{k(\fm)}\otimes_{k(\fm)}k(\fn)\cong
\aq n{R_{\fm}}{Q_{\fn\cap Q}}{k(\fm)}\otimes_{k(\fm)}k(\fn)
\]
which show that (i) implies (ii).  If $R\to S$ is the identity map, then
(ii) implies (i) by \cite[Suppl.\ Prop.\ 29]{hca}; the {\em proof\/} of
that proposition applies {\em verbatim\/} to every homomorphism $R\to S$
in which the ring $S$ is noetherian.
 \end{chunk}

We recall from \cite{lci} the definition of
l.c.i.\ homomorphism of noetherian rings.

\begin{definition}
\label{def: lci}
A homomorphism of noetherian rings $\vf\col R\to S$ is called {\em locally
complete intersection\/} (or {\em l.c.i.\/}) if for every $\fn\in\Spec
S$ the induced local homomorphism $\vf_\fn\col R_{\fn\cap R}\to S_\fn$
is complete intersection at $\fn S_\fn$ in the sense of \ref{def: ci}.
 \end{definition}

We return to the discussion, started in the Introduction, of criteria for
l.c.i.\ homomorphisms in terms of Andr\'e-Quillen homology.  Recall that
$\vf\col R\to S$ is said to be {\em locally of finite flat dimension\/}
if $\fd_R(S_\fn)<\infty$ for every $\fn\in\Spec S$.

\begin{chunk}
\label{flat lci}
If $\vf\col R\to S$ is a homomorphism of noetherian rings, then
the following conditions are equivalent.
\begin{enumerate}[{\quad\rm(i)}]
\item
$\aqdim SR < \infty$, and $\vf$ is locally of finite flat
dimension.
\item
$\aqdim SR\leq 1$\,.
\item
$\aq 2SR-=0$\,.
\item
$\aq nSR-=0$ for some $n\ge2$ with $\lfloor \frac n2\rfloor!$
invertible in $S$, and $\vf$ is locally of finite flat dimension.
\item
$\vf$ is locally complete intersection\,.
\end{enumerate}
Indeed, (v) is equivalent to (ii), (iii) by \cite[(1.2)]{lci},
and to (i), (iv) by \cite[(1.5)]{lci}.
\end{chunk}

A major difficulty in dealing with Quillen's Conjecture is that
$R$-algebras $S$ with $\aqdim SR\leq 2$ have not been described in
structural terms.  All known algebras satisfying this condition are
constructed by factoring some l.c.i.\ homomorphism.

\begin{example}
Let $Q\to S$ be an l.c.i.\ homomorphism from a noetherian ring $Q$ (say,
a surjection with kernel generated by a regular sequence $\bsg$), through
some homomorphism $Q\to R$ that is c.i.\ at all primes of $R$ contracted
from $S$ (say, a surjection with kernel generated by a regular sequence
in $(\bsg)$).  The desired vanishing property follows from \ref{flat lci},
via the Jacobi-Zariski exact sequence \cite[(5.1)]{hca}.
 \end{example}

The construction described above is rather rigid.  This is demonstrated
by the next result, which could be compared with another factorization
theorem for local homomorphisms from \cite{lci}:  If $\vf\circ\psi$
is l.c.i.\ and $\vf$ is locally of finite flat dimension, then $\vf$
is l.c.i.\ and $\psi$ is complete intersection at $\fn\cap R$ for each
$\fn\in \Spec S$.

\begin{theorem}
\label{lci factor}
If $\psi\col Q\to R$ and $\vf\col R\to S$ are homomorphisms of
noetherian rings such that $\vf\circ\psi$ is l.c.i., then the following
conditions are equivalent.
\begin{enumerate}[{\quad\rm(i)}]
\item
$\aqdim SR<\infty$.
\item
$\aqdim SR\leq 2$\,.
\item
$\aq 3SR-=0$\,.
\item
$\aq nSR-=0$ for some $n\ge 3$ such that $\lfloor \frac{n-1}2\rfloor!$
is invertible in $S$\,.
\item
$\psi$ is complete intersection at $\fn\cap R$ for each $\fn\in \Spec S$.
\end{enumerate}
 \end{theorem}

\begin{proof}
The Jacobi-Zariski exact sequence of Andr\'e-Quillen homology,
cf.\ \cite[(5.1)]{hca}, yields an exact sequence of functors on the
category of $S$-modules
\[
\cdots\lra\aq{n+1}SQ-\lra\aq{n+1}SR-\lra\aq{n}RQ-\lra\aq{n}SQ-\lra\cdots
\]
Since $\vf\circ\psi$ is l.c.i., we have $\aq{n}SQ-=0$ for $n\ge2$ by
\ref{flat lci}, so for $n\ge2$ we get isomorphisms $\aq{n+1}SR-\cong
\aq nRQ-$ of functors of $S$-modules.  Thus, each one of conditions (i),
(ii), and (iii) is equivalent to its primed version below:
\begin{enumerate}[{\quad\rm(i)}]
\item[\rm(i$'$)]
$\aq nRQ-=0$ on the category of $S$-modules for all $n\gg 0$.
\item[\rm(ii$'$)]
$\aq nRQ-=0$ on the category of $S$-modules for all $n\ge2$.
\item[\rm(iii$'$)]
$\aq 2RQ-=0$ on the category of $S$-modules.
\item[\rm(iv$'$)]
$\aq nRQ-=0$ on the category of $S$-modules for some $n\ge 2$ such that
$\lfloor \frac n2 \rfloor!$ is invertible in $S$.
\end{enumerate}

The equivalence of (ii$'$) and (v) results from \ref{localization} and
\ref{flat ci}.  For each $\fn\in\Spec S$, the homomorphism $\psi_{\fn\cap
R}\col Q_{\fn\cap Q}\to R_{\fn\cap R}$ is almost small by Corollary
\ref{small examples}.  The equivalence of (i$'$), (iii$'$), (iv$'$), and
(v) now comes from \ref{localization} and Theorem \ref{almost quillen}.
 \end{proof}

The last result implies part of \cite[(1.5)]{lci}: Quillen's Conjecture
holds when the ring $S$ is l.c.i., because the theorem then applies
with $Q=\BZ$.  More to the point, it reduces the proof of Theorem I from
the Introduction to a mere formality.

\begin{proof}[Proof of Theorem {\em I}]
Let $\psi\col S\to R$ be any section of $\vf$.

The map $\vf\circ\psi=\id_S$ is obviously l.c.i., so Theorem \ref{lci
factor} and Proposition \ref{split ci}.a show that the conditions of
Theorem I are equivalent.
 \end{proof}

It remains to deduce Theorem II from its local version established in
Section \ref{Andre-Quillen homology}.

\begin{proof}[Proof of Theorem {\em II}]
The implication (ii) $\implies$ (i) is clear.  By \eqref{deg1} and
\eqref{deg2}, the $S$-modules $\aq 1RSS$ and $\aq 2RSS$ are finitely
generated.  Thus the implications (i) $\implies$ (iii) $\iff$ (iv)
follow from Theorem \ref{local second} via the isomorphisms
\begin{gather*}
\Tor{\bu}RSS{}_\fn\cong\Tor{\bu}{R_{\fn\cap R}}{S_\fn}{S_\fn}\,;\\
\aq\bu SRS_\fn\cong\aq\bu{S_\fn}{R_{\fn\cap R}}{S_\fn}\,,
\end{gather*}
respectively of graded $S_\fn$-algebras and graded $S_\fn$-modules,
cf.\ \cite[(4.59), (5.27)]{hca}.

It remains to prove (iii) $\implies$ (ii). Let $D_\bu$ denote the
graded $S$-module with $D_n=\aq nSRS$ for $n=1,2$, and $D_n=0$
otherwise. The isomorphisms \eqref{deg1} and \eqref{deg2} define a
surjection $\tau\col\Tor{\bu}RSS \to D_\bu$ of graded $S$-modules. Since
$D_\bu$ is projective, we may choose an $S$-linear map $\theta\col
D_\bu\to\Tor{\bu}RSS$ with $\tau\circ\theta=\id_{D_\bu}$.  It extends
to a homomorphism of graded $S$-algebras
\[
\vartheta\col{\ts\bigwedge}_{S}D_1\otimes_S\sym_{S}D_2\lra\Tor{\bu}RSS\,.
\]
Theorem \ref{local second} and the isomorphisms above show that
$\vartheta_\fn$ is bijective for every $\fn\in\Max S$, so $\vartheta$
is an isomorphism.   It induces an isomorphism
\[
\vartheta\otimes_ST\col
({\ts\bigwedge}_{S}D_1\otimes_S\sym_{S}D_2)\otimes_ST\xra{\ \cong\ }
\Tor{\bu}RSS\otimes_ST
\]
of graded $T$-algebras.  As each $S$-module $\Tor nRSS$ is projective,
the Universal Coefficients Theorem shows that for every $S$-algebra $T$
the homomorphism
\[
\Tor{\bu}RSS\otimes_ST\xra{\ \cong\ }\Tor{\bu}RS{T}
\]
of graded $T$-algebras given by the K\"unneth map is bijective.
Composing the isomorphisms above we obtain the desired isomorphism of
graded $T$-algebras.
 \end{proof}

We finish the paper by revisiting the announcement in \cite{tata}.

\begin{remark}
\label{errata}
In \cite[(2.6)]{tata}, which is an avatar of Theorem \ref{lci factor},
instead of condition (\ref{lci factor}.i) one finds the seemingly weaker
condition: $\aqdim {S_\fn}R<\infty$ for all $\fn\in\Spec S$.  We show that
this condition is in fact equivalent to those in Theorem \ref{lci factor}.

Indeed, under the hypotheses of Theorem \ref{lci factor} (which coincide
with those \cite[(2.6)]{tata}) the composition of the maps $\psi\col Q\to
R$ and $R\to S_\fn$ is an l.c.i.\ homomorphism.   Fixing a prime ideal
$\fn$ of $S$, from Theorem \ref{lci factor} we see that $\aqdim {S_\fn}R$
is finite if and only if $\psi$ is c.i.\ at $\fq\cap R$ for all $\fq\in
\Spec S$ with $\fq\subseteq\fn$.  Letting now $\fn$ range over $\Spec
S$, we conclude that $\aqdim {S_\fn}R$ is finite for all $\fn\in\Spec S$
if and only if $\psi$ is c.i.\ at $\fq\cap R$ for all $\fq\in \Spec S$.
This is condition (\ref{lci factor}.v).

We take this opportunity to make two minor corrections to \cite{tata}:
\begin{enumerate}[\rm(1)]
\item[--]
in \cite[(2.6.i)]{tata} `at each $\fp\in\Spec R$ with $\fp\supseteq
\Ker(\vf)$' should read `at $\fn\cap R$ for each $\fn\in\Spec S$';
\item[--]
in \cite[line above (4.1)]{tata} `in [13, (4.1.iii)]'
should read `in [13, Theorem II]'.
\end{enumerate}
We thank the editors of the present article for pointing these errors
to us.
 \end{remark}

\end{document}